\documentclass[11pt,reqno]{amsart}
\usepackage{amsmath,amsfonts,amssymb,amscd,amsthm,amsbsy,bbm, epsf,calc,  comment, caption, enumerate}
\usepackage{color}
\usepackage[usenames,dvipsnames]{xcolor}
\usepackage{datetime}
\usepackage{soul}

\usepackage{thmtools, thm-restate}
\usepackage{thm-patch, thm-kv, thm-autoref}
\usepackage[colorlinks=true, urlcolor=blue, linkcolor=blue, citecolor=black]{hyperref}

\usepackage[pdftex]{graphicx}
\usepackage{wrapfig}
\usepackage{geometry}
\usepackage{scrextend}
\usepackage{mathtools}

\numberwithin{equation}{section}
\newcounter{hours}\newcounter{minutes}

\theoremstyle{plain}
\declaretheorem[title=Theorem, parent=section]{theorem}
\declaretheorem[title=Lemma,sibling=theorem]{lemma}

\declaretheorem[title=Proposition,sibling=theorem]{proposition}
\declaretheorem[title=Definition,sibling=theorem]{definition}
\declaretheorem[title=Remark,sibling=theorem]{rem}

\declaretheorem[title=Assumption,sibling=theorem]{assumption}

\newtheorem*{prop*}{Proposition}

\makeatletter
\newcommand\RedeclareMathOperator{%
	\@ifstar{\def\rmo@s{m}\rmo@redeclare}{\def\rmo@s{o}\rmo@redeclare}%
}
\newcommand\rmo@redeclare[2]{%
	\begingroup \escapechar\m@ne\xdef\@gtempa{{\string#1}}\endgroup
	\expandafter\@ifundefined\@gtempa
	{\@latex@error{\noexpand#1undefined}\@ehc}%
	\relax
	\expandafter\rmo@declmathop\rmo@s{#1}{#2}}
\newcommand\rmo@declmathop[3]{%
	\DeclareRobustCommand{#2}{\qopname\newmcodes@#1{#3}}%
}
\@onlypreamble\RedeclareMathOperator
\makeatother

\def\ep{\varepsilon}
\def\al{\alpha}
\def\del{\delta}

\def\om{\omega}

\def\Om{\Omega}
\def\gam{\gamma} 
\def\Gam{\Gamma}
\def\lam{\lambda}
\def\Lam{\Lambda}

\def\grad{\nabla}
\def\integer{\mathbb Z}
\def\Natural{\mathbb N}
\def\real{\mathbb R}
\def\A{\mathcal A}
\def\B{\mathcal B}
\def\C{\mathcal C}

\def\K{\mathcal K}

\def\M{\mathcal M}

\def\P{\mathcal P}
\def\R{\mathbb R}

\def\S{\mathcal S}

\def\union{\cup}
\def\Union{\bigcup}

\def\Expectation{\mathbb{E}}
\def\Indicator{{\mathbbm{1}}}

\def\Prob{\mathbb{P}}

\def\Tr{\textnormal{tr}}

\RedeclareMathOperator{\div}{\textnormal{div}}

\def\polhk#1{\setbox0=\hbox{#1}{\ooalign{\hidewidth
			\lower1.5ex\hbox{`}\hidewidth\crcr\unhbox0}}}

\newcommand{\abs}[1]{\left| #1 \right|}

\newcommand{\norm}[1]{\lVert#1\rVert}

\newcommand{\sub}{\mathrm{sub}}
\newcommand{\super}{\mathrm{super}}

\begin{document}

	\title{Hamilton-Jacobi-Bellman equations on graphs}

	\author{Nicol\`o Forcillo}
	\author{Jun Kitagawa}
	\author{Russell W. Schwab}

	\address{Department of Mathematics\\
		Michigan State University\\
		619 Red Cedar Road \\
		East Lansing, MI 48824}
	\email{ forcill1@msu.edu, kitagawa@math.msu.edu, rschwab@math.msu.edu}

	\begin{abstract}

	Here, we study Hamilton-Jacobi-Bellman equations on graphs.  These are meant to be the analog of any of the following types of equations in the continuum setting of partial differential and nonlocal integro-differential equations: Hamilton-Jacobi (typically first order and local), Hamilton-Jacobi-Bellmann-Isaacs (first, second, or fractional order), and elliptic integro-differential equations (nonlocal equations).  We give conditions for the existence and uniqueness of solutions of these equations, and work through a long list of examples in which these assumptions are satisfied.  This work is meant to accomplish three goals: complement and unite earlier assumptions and arguments focused more on the Hamilton-Jacobi type structure; import ideas from nonlocal elliptic integro-differential equations; and argue that nearly all of the operators in this family enjoy a common structure of being a monotone function of the differences of the unknown, plus ``lower order'' terms.  This last goal is tied to the fact that most of the examples in this family can be proven to have a Bellman-Isaacs representation as a min-max of linear operators with a graph Laplacian structure.

	\end{abstract}

	\date{\today,\ arXiv version 1.}

	\thanks{
	J. Kitagawa was supported in part by National Science Foundation grant DMS-2246606. R. Schwab acknowledges support from the Simons Foundation for a Travel Support for Mathematicians grant.}
	
	
	\maketitle
	
	\markboth{Hamilton-Jacobi-Bellman equations on graphs}{Elliptic equations on graphs}



	
\section{Introduction and main result}\label{sec:Intro}

In this paper, we study Hamilton-Jacobi-Bellman equations on graphs, and present a class of them which exhibit unique solutions.  These equations are intended to be analogous to those in the setting of partial differential equations classified as either Hamilton-Jacobi equations, Bellman equations, or (nonlinear) elliptic equations.   
Methods from the PDE theory have already been used to study solutions to certain equations that have a strong resemblance to Hamilton-Jacobi equations (such as eikonal or Bellman type equations, see for example, \cite{Cald_Ette}, \cite{Ober-2006ConvDiffSchemesSIAM}).  Here we present a complementary point of view utilizing methods more along the lines of elliptic and integro-differential PDE, which we will explain below.  The main defining feature of our equations is that the corresponding operators enjoy something we will call the global comparison property, defined below.

The set-up is the following.  Let $G$ be a graph with $N$ vertices, which we denote by
\begin{align*}
	G=\{x_1,\ldots,x_N\}.
\end{align*}
For the moment, we do not add any edges on $G$, but we will consider edges with weights later as necessary. 
We work with the set of all real-valued functions on $G$ and which we denote by
	\begin{align}\label{eqIntro:CG}
		C(G)\coloneqq  \{u\ |\ u: G\to \R\}.
	\end{align}
	Note that one can view $C(G)\cong \R^N$ by identifying a function $u\in C(G)$ with the vector $(u(x_1),\ldots, u(x_N))\in\R^N$.
	We will let $I$ be a general operator acting on $C(G)$, i.e.
		\begin{align*}
			I: C(G)\to C(G),\ \ \ \text{and we use the notation}\ \ \ 
			I(u,x)\in\real,\ u\in C(G),\ x\in G.
		\end{align*}

	We are interested in studying equations on graphs of the form
	\begin{equation}\label{eqIntro:HJBGeneric}
		\begin{cases}
			I(u,\cdot)=f& \mbox{on }G\setminus \Gamma,\\
			u=g & \mbox{on } \Gamma,
		\end{cases}
	\end{equation}
	where  
	\begin{align*}
		\Gamma\subset G,\ \ \  
		\Gamma\neq \emptyset,\ \ \ 
		\Gamma\neq G, 
	\end{align*}
	with
	\begin{align*}
		f \in C(G\setminus \Gam)\ \ \ \text{and}\ \ \ 
		g\in C(\Gam).
	\end{align*}
	The functions $f$ and $g$ play the roles, respectively of source and boundary data.

	We will expand upon many examples of $I$ later in the paper, but for now, a canonical class of linear examples are $I$ of the form
	\begin{align}\label{eqIntro:IGraphLaplace}
		I(u,x) = \sum_{y\in G} a(x, y)(u(y)-u(x)),
	\end{align}
	where $a: G\times G\to \R$, $a\ge 0$.  A canonical class of nonlinear examples would be for some index sets, $\A$ and $\B$, and a doubly indexed class of linear operates,
	\begin{align}\label{eqIntro:MinMax}
		I(u,x) = \min_{\al\in\A}\max_{\beta\in \B}
		\left( f^{\al\beta}(x) + 
		c^{\al\beta}(x)u(x) + \sum_{y\in G} a^{\al\beta}(x, y)(u(y)-u(x))
		\right).
	\end{align}
	We will expand upon this in Section \ref{sec:Examples}.  It turns out that in most circumstances, the structure in \eqref{eqIntro:MinMax} is generic, in the sense that most elliptic operators that satisfy our assumptions for solving the equation \eqref{eqIntro:HJBGeneric} must in fact admit a min-max representation in \eqref{eqIntro:MinMax} (at least those operators that are locally Lipschitz on $C(G)$).   This will be discussed in subsections \ref{subsec:GeneralMinMax} and \ref{subsec:GeneralH-CE}.   Basically, operators in \eqref{eqIntro:IGraphLaplace} and \eqref{eqIntro:MinMax} are the direct analogs of linear and fully nonlinear (degenerate) elliptic nonlocal integro-differential operators, and we use the analogy to present a theory for Hamilton-Jacobi-Bellman equations on graphs that follows along the lines of the elliptic nonlocal equations in the continuum setting.

 Our investigations of $I$, are based on a feature we call the \emph{global comparison property}.
	
	\begin{definition}\label{defIntro:GCP}
		We say that $I:C(G)\to C(G)$ satisfies the global comparison property (GCP), provided that for $u,v\in C(G)$,
		\begin{align*}
			\text{whenever}\ u\leq v\text{ on }G\ \text{and}\ u(x_0)=v(x_0),\ \ 
			\text{then}\ \ 
			I(u,x_0)\leq I(v,x_0).
		\end{align*}
	\end{definition}

	\begin{rem}
		
		In case $I$ is of the form \eqref{eqIntro:IGraphLaplace}, respectively \eqref{eqIntro:MinMax}, nonnegativity of $a$, respectively $a^{\al\beta}$, is equivalent to $I$ having the GCP.
		 
	\end{rem}

	We are interested in identifying situations in which the equation \eqref{eqIntro:HJBGeneric} admits a unique solution-- i.e. assumptions on $I$, $f$, and $g$.  Generically, it is typical to break this into two questions: when solutions exist, are they unique; and when do solutions exist?  Historically, uniqueness (at least in the context of Hamilton-Jacobi equations) is often resolved via a result showing that a subsolution and a supersolution that are ordered on the boundary of the domain must also be ordered inside the domain, hence is typically called a \emph{comparison result}. 
	
	To this end, we give a few definitions. 
    \begin{definition}[Sub/supersolution]\label{defIntro:Sub-Super-Solution}
    	
    	We say that $u$ is a \textit{subsolution, (respectively strict subsolution)} to $I(u,\cdot)=f$ on $G\setminus\Gam$, provided 
		\begin{align*}
			I(u,\cdot)\geq f\ \text{on}\ G\setminus \Gam,\ \ \text{(respectively}\ I(u,\cdot)> f\text{on}\ G\setminus \Gam).
		\end{align*}

		Similarly, we say that $v$ is a \textit{supersolution, (respectively strict supersolution)} to $I(u,\cdot)=f$ on $G\setminus\Gam$, provided 
		\begin{align*}
			I(v,\cdot)\leq f\ \text{on}\ G\setminus \Gam,\ \ \text{(respectively}\ I(u,\cdot)< f\text{on}\ G\setminus \Gam).
		\end{align*}

    \end{definition}

	\begin{rem}
		We point out the very important fact that the sign convention taken here differs from that in many works, e.g. \cite{Cald_Ette}, with regards to a sub and supersolution.  Often, the GCP is assumed to have the reverse inequality, and consequently a subsolution would be a function with $I(u, \cdot)\leq f$.  In the continuous case, this corresponds to the difference between considering $I(u, \cdot)=\Delta u$ verses $I(u, \cdot)=-\Delta u$ (as an example, note the differing sign conventions in \cite{CaCa-95} and \cite{CrandalIshiLions-92UsersGuide}).  
	\end{rem}

As a technical point, we will use the $L^\infty$ norm to describe convergence in $C(G)$.  Since we can associate $C(G)$ with $\real^N$, all norms are equivalent, but we choose one for convenience.

\begin{definition}
	We say that for $u_k,u\in C(G)$, $u_k\to u$, provided
	\begin{align*}
		\lim_{k\to\infty}\norm{u_k-u}_{L^\infty(G)}=0,
	\end{align*}
	where
	\begin{align*}
		\norm{w}_{L^\infty(G)} = \max_{y\in G}\abs{w(y)},
	\end{align*}
	and we will subsequently denote $\norm{w}_{C(G)} = \norm{w}_{L^\infty(G)}$. 
\end{definition}

	Our goal is to give a reasonable list of assumptions (detailed in Subsection~\ref{subsec:Assumptions} below) under which the following uniqueness and existence results will be true for \eqref{eqIntro:HJBGeneric}.

	\begin{theorem}[Comparison/uniqueness result]\label{thmIntro:Comparison}
		
		Suppose $I$ satisfies Assumption \ref{assume:Uniqueness} (below) and $u$, $v\in C(G)$ satisfy
		\begin{align*}
			I(u, \cdot)\geq  I(v, \cdot) \ \ \text{on}\ G\setminus\Gam.
		\end{align*}
		Additionally, suppose one of the two following conditions holds: 
		 
			\begin{enumerate}[(a)]
				\item\label{assume:StrictSubSolPerturbation} (Strict subsolution perturbation) There exists a sequence $\{u_k\}_{k=1}^\infty\subset C(G)$,  such that for all $k\in \mathbb{N}$,
			\begin{align*}
				I(u_k,\cdot)>I(u,\cdot)\quad\mbox {on }G\setminus\Gamma,
			\end{align*}
			and
			\begin{align*}
				\lim_{k\to\infty} \norm{u-u_k}_{C(G)}=0.
			\end{align*}

				\item\label{assume:StrictSuperSolPerturbation} (Strict supersolution perturbation) There exists a sequence $\{v_k\}_{k=1}^\infty\subset C(G)$,  such that for all $k\in \mathbb{N}$,
			\begin{align*}
				I(v_k,\cdot)<I(v,\cdot)\quad\mbox {on }G\setminus\Gamma,
			\end{align*}
			and
			\begin{align*}
				\lim_{k\to\infty} \norm{v-v_k}_{C(G)}=0.
			\end{align*}
			\end{enumerate}
		 Then 
		\begin{align*}
			\max_{G} (u-v)_+\le \max_{\Gamma} (u-v)_+,
		\end{align*}
		where $(u-v)_+=\max\{u-v,0\}$.
		 	
	\end{theorem}
	
	We note that this is often interpreted in terms of subsolution and supersolution to a particular equation, in the sense that for some $f\in C(G)$,
	\begin{align*}
		I(u, \cdot)\geq f\ \ \text{and}\ \ \ I(v, \cdot)\leq f\ \ \text{on}\ G\setminus\Gam,
	\end{align*}
	and then Theorem \ref{thmIntro:Comparison} says that the maximum of $u-v$ will occur on the ``boundary'' of the domain.

	\begin{theorem}[Existence]\label{thmIntro:Existence}
		
		Suppose $I$ satisfies Assumptions \ref{assume:Uniqueness} and \ref{assume:Existence} (below), and $f$, $\Gamma$, and $g$ as in \eqref{eqIntro:HJBGeneric} are such that there exist $U_{f,g,\sub}$, $V_{f,g,\super}\in C(G)$ satisfying
			\begin{align}
				\begin{split}\label{eqn: special sub/super}
					I(U_{f,g,\sub},\cdot)&\geq f \ \text{on}\ G\setminus\Gam
					\ \text{and}\ \ U_{f,g,\sub}\leq g \ \ \text{on}\ \Gam,\\
					I(V_{f,g,\super},\cdot)&\leq f \ \text{on}\ G\setminus\Gam
					\ \text{and}\ \ V_{f,g,\super}\geq g \ \ \text{on}\ \Gam.
				\end{split}
			\end{align}
		
		If $I$ also satisfies the conclusion of Theorem \ref{thmIntro:Comparison}, i.e.
		\begin{align}\label{eqIntro:ExistenceTheorem-ComparisonAssumption}
			&I(u,\cdot)\geq I(v,\cdot)\ \text{in}\ G\setminus\Gam\   
			\text{and}\ u\leq v\ \text{on}\ \Gam\\
			&\implies\ u\leq v\in G,\nonumber
		\end{align}
		then there exists $u\in C(G)$, where $u$ solves \eqref{eqIntro:HJBGeneric}.
		
	\end{theorem}

\begin{rem}\label{rem: alternate existence condition}
	We note that to prove the existence result Theorem \ref{thmIntro:Existence}, we assume the comparison result Theorem \ref{thmIntro:Comparison} holds between \emph{all} pairs of sub and supersolutions, which is typical of the presentation in the continuum case of partial differential and integro-differential equations. However, we note it is possible to obtain the existence result by assuming some conditions on a more restricted set of sub or supersolutions. In particular, one can replace \eqref{eqIntro:ExistenceTheorem-ComparisonAssumption} and still obtain existence: 
	
			\begin{addmargin}{1cm}
			If either $V_{f,g,\super}$ satisfies condition \eqref{assume:StrictSuperSolPerturbation} from Theorem \ref{thmIntro:Comparison}, or every subsolution $u\in C(G)$ such that 
			\begin{align*}
					I(u,\cdot)\geq f \ \text{on}\ G\setminus\Gam
					\ \text{and}\ \ u\leq g \ \ \text{on}\ \Gam
				\end{align*}
	satisfies condition \eqref{assume:StrictSubSolPerturbation} from Theorem \ref{thmIntro:Comparison},
			then there exists a unique solution to \eqref{eqIntro:HJBGeneric}.
			\end{addmargin}
	
\end{rem}

\begin{rem}[Elliptic]
	
	Above, we mentioned that the focus of this work is to study elliptic equations on graphs and to find when there exist unique solutions.  In the PDE setting, it is much more straightforward to determine what should be an elliptic equation.  In this context, one could say that that combination of the operator enjoying the GCP, and its related equation enjoying the outcome of Theorem \ref{thmIntro:Comparison} is what we mean by an elliptic equation for functions on the graph, $G$.  This is, of course, not precise, but it tracks closely with the PDE theory in the continuum setting, where our equations would be referred to as \emph{proper} or \emph{degenerate elliptic}, as in \cite{CrandalIshiLions-92UsersGuide}.
	
\end{rem}


\subsection{Related results and literature}\label{subsec:RelatedResultsLiterature}

\subsubsection{Relationship to \cite{Cald_Ette}}

Our motivation for this work comes from \cite{Cald_Ette}, and we hope that this serves as a complementary approach to Hamilton-Jacobi equations on graphs that was presented therein.  We will discuss two important examples from \cite{Cald_Ette} and \cite{BungertCalderRoith-2023UnifRatesLipLearningOnGraphs} later, in subsection \ref{subsec:Eikonal}, and then more generally in subsection \ref{subsec:GeneralH-CE}, but we want to make some initial comments here.  In \cite{Cald_Ette}, the approach to identifying equations like \eqref{eqIntro:HJBGeneric} that should be considered versions of Hamilton-Jacobi equations on graphs came from the fact that in the continuum setting, for $u:\real^d\to\real$, Hamilton-Jacobi equations typically have the structure 
\begin{align*}
	H(\grad u(x), u(x), x ) = f(x)\ \ \text{in some open set,}\ \ D\subset\real^d.
\end{align*}
Notably, these equations are local and first order, among other things.
In the first order case, since for any two differentiable functions, say $u$ and $v$ with $u\leq v$ and $u(x_0)=v(x_0)$, it follows that $\grad u(x_0)=\grad v(x_0)$, the only monotonicity in the equation that needs to be assumed is with respect to the pointwise value $u(x)$, and typically $H$ is increasing in $u$ (this is referred to as the equation is ``proper'' or ``degenerate elliptic'' as in \cite{CrandalIshiLions-92UsersGuide}).  In the continuum setting, \emph{for first order equations}, this observation is why the GCP does not play a substantial role in the theory (as in stark contrast to the second order elliptic PDE theory).  However, as one can see from the proof of Theorem \ref{thmIntro:Comparison}, some form of monotonicity property of the equation \eqref{eqIntro:HJBGeneric} is needed for the theory related to the Hamilton-Jacobi and elliptic equations when working in the discrete setting of a graph, and this becomes an assumption on $H$. 
In \cite{Cald_Ette}, the authors define a natural candidate for the ``gradient'' of a function, $u\in C(G)$, which at $x_i$ is given by the vector of differences,
\begin{align*}
	\grad_G u(x_i) = \left(  
	u(x_i)-u(x_1),
	u(x_i)-u(x_2),
	\dots,
	u(x_i)-u(x_N)
	\right),
\end{align*}
and then monotonicity of the equation in this $\grad_G u$ variable is assumed.  That is to say, they introduce a Hamilton-Jacobi equation on a graph as some function $H$, with
\begin{align*}
	H:\real^N\times\real\times G \to\real,
\end{align*}
and
\begin{align}\label{eqIntro:CalderHJ-Eq}
	H(\grad_G u(x), u(x), x) = f(x)\ \ \text{on}\ \ G\setminus\Gam,
\end{align}
where $H$ is monotone in its first two arguments; this monotonicity is discussed later, in subsection \ref{subsec:GeneralH-CE}.  
We also note, this idea of considering a discrete Hamilton-Jacobi equation by replacing the gradient with a vector of differences was already present in \cite{Ober-2006ConvDiffSchemesSIAM}, and other numerical stencil based schemes. It turns out that when viewed in terms of the GCP and elliptic operators, this choice of the structure like \eqref{eqIntro:CalderHJ-Eq}, is not just natural, but it is more or less required of any reasonable operator which has the GCP (which, one naturally expects from discrete Hamilton-Jacobi equations).  In fact, the family of equations for which a comparison theorem will hold effectively requires $H$, or $I$ in our context here, to be a monotone function of the differences of $u$.  This stems from the work in \cite{GuSc-2019MinMaxNonlocalCALCVARPDE} and \cite{GuSc-2019MinMaxEuclideanNATMA}, is explained in more detail in subsection \ref{subsec:GeneralMinMax}, below.

\subsubsection{Other works}

In \cite{BungertCalderRoith-2023UnifRatesLipLearningOnGraphs} and \cite{Cald_Ette}, the notion of studying Hamilton-Jacobi equations on a graph arose naturally to try to compute certain properties of graphs, like distances and other related objects.  In this context, the graph was the original object of study, and theory from partial differential equations in the continuum setting was imported to guide the relevant equations on a graph in their setting.

There is a related reason to study equations on certain graphs for which the flow of information goes in the other direction, namely discretization schemes for numerically approximating solutions to partial differential equations.  For elliptic equations this goes all the way back to \cite{MotzkinWasow-1952ApproxLinDiffEq-JMP}, and there are natural notions of discretization methods which then lend themselves to solving a related equation in the discretized space.  So in this context, a PDE in the continuum setting is the original object, and the equation on the graph is meant to be a method for approximating solutions.  In the linear and elliptic setting, some very detailed properties of solutions of discrete elliptic equations were given in \cite{KuoTrudinger-1990-LinearDifferenceInequalitiesRandomCoefficient-MathComp} and \cite{KuoTrudinger-1992DiscreteMethodsNonlinearSIAMNumAnal}.  For degenerate elliptic and parabolic, possibly nonlinear equations, and first order Hamilton-Jacobi equations, in the continuum setting, one effective choice for non-classical solutions is that of \emph{viscosity solutions}, for which some standard introductory references are \cite{BardiCapuzzoDolcetta-1997-ViscSolBookOptControlBellman-Birkhauser}, \cite{CrandalIshiLions-92UsersGuide}, \cite{ImbertSilvestre-2013IntroToFullyNonlinearParabolic}.  One of the many reasons that viscosity solutions are widely used is that they lend themselves very well to certain monotone approximation schemes.  We note that monotonicity of approximation schemes is intimately connected to the GCP we use herein (Definition \ref{defIntro:GCP}), as well as the fact that these equations respect the comparison of sub and super solutions.  A fundamental result in this direction that gave way to a large class of solutions of equations on a lattice, and their limiting solutions, is the convergence of schemes to viscosity solutions shown in \cite{BarlesSouganidis-1991ConvApproxSchemesFullyNonlinear}.  A good reference of earlier work on Hamilton-Jacobi equations where one can see this notion of solving for viscosity solutions on graphs is \cite{Soug-85}.  One example of an influential breakthrough of viscosity solutions to other previously unrelated equations is in the context of mean curvature flow, for which some good references are \cite{ChenGigaGoto-1991ViscSolMeanCuravture-JDG}, \cite{EvansSpruck-1991-MCM-JDG}, \cite{EvansSonerSouganidis-1992-PhaseTransistionAndMCM-CPAM}.  This led to more detailed exploration for these approximation schemes in the context of mean curvature flow in, e.g. \cite{CrandallLions-1996ConvDiffSchemeNonlineParabolicMCFlow-NumerischeMath}.  This theory for schemes on graphs was distilled more carefully into the notion of viscosity solutions on graphs in the work of \cite{Ober-2006ConvDiffSchemesSIAM}.  Furthermore, the idea of using a replacement for mean curvature flow on graphs appeared in \cite{VanGennipGuillenOstingBertozzi-2014-MCThresholdMBOOnGraphs-MilanMath}, and this can be considered closely aligned to the ideas we present here.

Since one of our goals is to import ideas from nonlocal integro-differential equations to the graph setting, we will mention some references related to nonlocal integro-differential equations.  These equations are related to the generators and control/game problems for pure jump Markov processes.
A basic example of a linear elliptic integro-differential equation for a domain, $D\subset\real^d$, and unknown, $u:\real^d\to\real$, is described by a drift, $b:\real^d\to\real^d$, a nonlocal difference interaction, $K:\real^d\times\real^d\to[0,\infty)$, and complementary data, $g:\real^d\setminus D\to\real$, with the form
\begin{align}\label{eqIntro:LinearIntDiff}
	\begin{cases}
		b(x)\cdot \grad u(x) + \int_{\real^d} (u(x+h)-u(x) - \Indicator_{B_1}(h)\grad u(x)\cdot h) K(x,h)dh
		=f(x)\ &\text{in}\ D,\\
		u(x) = g(x)\ &\text{on}\ \real^d\setminus D.
	\end{cases}
\end{align}
In this context, $K(x,h)dh$ would be called the L\'evy measure for the operator, and at a most basic level, it satisfies
\begin{align*}
	\sup_x \int_{\real^d}\min\{\abs{h}^2,1\}K(x,h)dh <+\infty,
\end{align*}
from which one can see, via Taylor expansion, that the integro-differential operator is well defined when $u\in C^{1,1}(D)$.
Here, $K(x,\cdot)$ describes the interaction of differences of $u$, based at $x$, and so a translation invariant version of the linear equation would have $b$ is a constant and $K(x,h)$ is independent of $x$.  One of the simplest canonical choices of $K$ is, for $\al\in(0,2)$,
\begin{align*}
	K(h) = c_{d,\al}\abs{h}^{-d-\al},
\end{align*}
where $c_{d,\al}$ is a dimensional constant proportional to $(2-\al)$, and in this case, when $b\equiv 0$, the equation becomes
\begin{align*}
	-(-\Delta)^{\al/2}u(x) = f(x)\ \ \text{in}\ \ D,
\end{align*}
where $(-\Delta)^{\al/2}$ is the $\al/2$ power of the non-negative operator, $(-\Delta)$ (see, e.g. \cite{Land-72}).
Nonlinear elliptic nonlocal equations are essential aspects of the Hamilton-Jacobi-Bellman theory, and generally they take a min-max form similar to (see, e.g. \cite{CaSi-09RegularityIntegroDiff}, \cite{GuSc-2019MinMaxNonlocalCALCVARPDE}, \cite{GuSc-2019MinMaxEuclideanNATMA})
\begin{align}\label{eqIntro:IntDiffNonlinear}
	\min_{i}\max_{j} \left(
	b^{ij}(x)\cdot \grad u(x) + \int_{\real^d} (u(x+h)-u(x) - \Indicator_{B_1}(h)\grad u(x)\cdot h) K^{ij}(x,h)dh
	\right)
	=f(x)\ \ \text{in}\ \ D,
\end{align}
where $b^{ij}$ and $K^{ij}$ are a family of coefficients, each of which is as in (\ref{eqIntro:LinearIntDiff}).  Some early references for this family of equations includes: \cite{Bass1988-UniquenessInLawPureJumpMarkov-ProbThRelField}, \cite{BaKa-05Holder}, \cite{BaLe-2002Harnack}, \cite{Kass-2009APrioriIntDiff}, \cite{Silv2006-HolderIntDiff-IUMJ}.  In fact, there is even a work that studies the properties of solutions to discrete versions of (\ref{eqIntro:LinearIntDiff}) on $\integer^d$ (albeit in the setting of divergence form operators) in \cite{BaLe-2002TransitionProb}, where heat kernel estimates are obtained via probabilistic methods that analyze related stochastic processes on $\integer^d$.

One of the reasons, among many, that equations  like (\ref{eqIntro:LinearIntDiff}) and (\ref{eqIntro:IntDiffNonlinear}) have received so much attention in the past 25 or so years is that when $K$ is properly weighted with a factor of $(2-\al)$, like the $\al/2$ Laplacian, above, then all second order elliptic linear and nonlinear equations can be captured as limiting cases of (\ref{eqIntro:LinearIntDiff}) and (\ref{eqIntro:IntDiffNonlinear}), in the same way that $(-\Delta)^{\al/2}$ converges to $(-\Delta)$ as $\al\to 2$.   That is to say, in the linear setting, one captures as limiting cases, equations of the form
\begin{align*}
	b(x)\cdot\grad u(x) + \Tr(A(x)D^2u(x)) = f(x)\ \ \text{in}\ \ D,
\end{align*} 
where $A\geq 0$ is a symmetric non-negative matrix.
An example of the power of this robustness in studying the family of equations like (\ref{eqIntro:IntDiffNonlinear}) can be seen in, e.g. \cite{CaSi-09RegularityIntegroDiff}, \cite{ChangLaraDavila-2016HolderNonlocalParabolicDriftJDE}, \cite{SchwabSilvestre-2014RegularityIntDiffVeryIrregKernelsAPDE}.  Equations like (\ref{eqIntro:LinearIntDiff}) also have analogs in the structure of divergence form equations, and some robust estimates in that context are \cite{FelsingerKassmann-2013LocalRegForNonlocalEqCPDE}, \cite{KassSchwa-2013RegularityNonlocalParaParma}.

One of the main points of this paper is that in the setting for which the GCP is applicable, most Hamilton-Jacobi-Bellman equations on $G$, like (\ref{eqIntro:HJBGeneric}), can be treated as discrete versions, with similar methods, of equations like (\ref{eqIntro:LinearIntDiff}) and (\ref{eqIntro:IntDiffNonlinear}).


\subsection{Assumptions}\label{subsec:Assumptions}

	We now state the assumptions that will lead to Theorems \ref{thmIntro:Comparison} and \ref{thmIntro:Existence}.  Additionally, we will address how Conditions \eqref{assume:StrictSubSolPerturbation} or \eqref{assume:StrictSuperSolPerturbation} from Theorem~\ref{thmIntro:Comparison} can be satisfied in practice, for a wide class of operators $I$.

We can now state the assumptions.

	\begin{assumption}[For comparison/uniqueness]\label{assume:Uniqueness} 
		
		The assumptions on $I$ related to uniqueness of solutions are as follows.
		
		\begin{enumerate}[(i)]

			\item\label{assume:GCP}
			
			$I$ has the GCP, given in Definition \ref{defIntro:GCP}
			
			\item\label{assume:SubtractConstant}
			
			For  all $u\in C(G)$, and any constant function $c\geq 0$ on $G$,
			\begin{align*}
				I(u-c,\cdot)\ge I(u,\cdot).
			\end{align*}	
			
		\end{enumerate}
		
	\end{assumption}

	\begin{assumption}[For existence]\label{assume:Existence}
		
		We require the following conditions for the functional $I$:
		
		\begin{enumerate}[(i)]

			\item\label{assume:UpperSemicon} For each $x\in G$, fixed,  $I(\cdot,x)$ is upper semicontinuous (u.s.c.) as a functional, $I(\cdot,x): C(G)\to\real$.

			\item\label{assume:LSCBump} (Positive perturbation property)  For each $x_0\in G$, fixed,  defining a particular bump function, $b_{x_0}$ as 
			
			\begin{align*}
				b_{x_0}(x)\coloneqq 
				\begin{cases}
					1 & \mbox{if }x=x_0,\\
					0 & \mbox{if }x\neq x_0,
				\end{cases}
			\end{align*}
			we require 
			\begin{align*}
				\lim_{t\to 0^+}\inf_{x \in G\setminus \Gamma}(I(u+tb_{x_0},x)-I(u,x))\ge 0,\quad\mbox{for any }u\in C(G).
			\end{align*}
			
		\end{enumerate}
		
	\end{assumption}


\subsection{Comments on the assumptions}\label{subsec:CommentsOnAssumptions}

As should become evident from the eventual proof of Theorem \ref{thmIntro:Comparison}, it is very easy to show that a strict subsolution and a (not necessarily strict) supersolution will be ordered in the domain if they are ordered on the ``boundary'' set $\Gam$.  In the continuous setting, often the method used for first order local Hamilton-Jacobi equations differ from those used for (degenerate) elliptic and parabolic local and nonlocal equations.  All of those equations enjoy the GCP, and they are amenable to viscosity solutions techniques, and so those arguments are usually similar, but the ways in which the various assumptions are satisfied are often different.  One feature we hope to illustrate in this article is that in the discrete graph setting, the two usually different types of equations become much more similar.

First, we will list a number of ways in which Conditions \eqref{assume:StrictSubSolPerturbation} or \eqref{assume:StrictSuperSolPerturbation} from Theorem~\ref{thmIntro:Comparison} can be satisfied. Most of the hypotheses below are taken from \cite{Cald_Ette}.

\begin{lemma}\label{lemIntro:StrictSubsolPerturbation}
	 
 	If $u\in C(G)$ and $I$ satisfy any of the following conditions, then $u$ satisfies Condition \eqref{assume:StrictSubSolPerturbation} from Theorem~\ref{thmIntro:Comparison}:
	\begin{enumerate}[(i)]
		\item \label{assumptionLemIntro:StrictProper} For any constant $c>0$, $I(u-c,\cdot)>I(u,\cdot)$ on $G$. 
		\item \label{assumptionLemIntro:PositiveHomogeneous}For some $p>0$, $I$ is positively $p$-homogeneous in $u$ and $I(u, \cdot)>0$.
		  \item \label{assumptionLemIntro:BarrierFunction}There exists $\phi_u\in C(G)$ and a sequence $\epsilon_k\searrow 0$ such that,
		\begin{align*}
			I(u+\epsilon_k\phi_u, \cdot)> I(u, \cdot)\ \ \text{on}\ \ G\setminus \Gam.
		\end{align*}
		\item \label{assumptionLemItro:ConcavePlusSubsolution}For each $x\in G$, the function $I(\cdot, x)$ is concave on $C(G)$, and there exists $\phi_u\in C(G)$ and a constant $\lam_u>0$, with
		\begin{align*}
			I(\phi_u, \cdot) \geq \lam_u + I(u, \cdot).
		\end{align*} 
		
	\end{enumerate}
\end{lemma}

\begin{rem}
	We emphasize here that the ``barrier function'', $\phi_u$ in \eqref{assumptionLemIntro:PositiveHomogeneous} above is allowed to depend upon $u$.  This is in contrast with the continuum setting of elliptic PDE, where typically one is fortunate enough to find a fixed function $\phi$ depending only on the domain and the operator, $I$ (such as a quadratic or exponential function in many cases).
\end{rem}

\begin{proof}[Proof of Lemma \ref{lemIntro:StrictSubsolPerturbation}]
	If $I$ satisfies \eqref{assumptionLemIntro:StrictProper}, then we can define $u_k = u -\frac{1}{k}$.
	If $u$ satisfies Condition \eqref{assumptionLemIntro:PositiveHomogeneous}, simply take $u_k\coloneqq u\cdot(k+1)/k$, and if $u$ satisfies Condition \eqref{assumptionLemIntro:BarrierFunction} let $u_k\coloneqq u+\epsilon_k \phi_u$.
	
	Now suppose $u$ satisfies Condition \eqref{assumptionLemItro:ConcavePlusSubsolution}, then for each $x\in G$,
	\begin{align*}
		I(k^{-1} u + (1-k^{-1})\phi_u, x)\geq k^{-1} I(u, x) + (1-k^{-1})I(\phi_u, x)
		\geq \lam_u + I(u, \cdot)>I(u, \cdot),
	\end{align*}
hence we may take $u_k\coloneqq k^{-1} u + (1-k^{-1})\phi_u$. finishing the proof.	
\end{proof}
\begin{rem}\label{remIntro:ConvexSuperSolutionPerturbation}
 Under appropriate conditions mirroring those of Lemma \ref{lemIntro:StrictSubsolPerturbation} above, one can show when a function will satisfy \eqref{assume:StrictSuperSolPerturbation} from Theorem~\ref{thmIntro:Comparison}. For example, if for each $x\in G$, the function $I(\cdot, x)$ is convex on $C(G)$, and there exists $\phi_v\in C(G)$ and a constant $\lam_v>0$, with
		\begin{align*}
			I(\phi_v, \cdot) \leq -\lam_v + I(v, \cdot).
		\end{align*} 
		then $v$ satisfies \eqref{assume:StrictSuperSolPerturbation} from Theorem~\ref{thmIntro:Comparison}.
\end{rem}

A convenient feature of sub/supersolutions is that their class is preserved under taking finite maxima/minima.

    \begin{lemma}\label{lemma:MaxSubsol}
    	If $u_1$ and $u_2$ are two subsolutions to $I(u,\cdot)=f$ on $G$, then $w\coloneqq  \max\{u_1,u_2\}$ is a subsolution to $I(u,\cdot)=f$ on $G$.
    \end{lemma}

    \begin{proof}
    	
    	The proof is a direct consequence of GCP. Indeed, by definition of $w$, for any $x\in G$, we have $w(x)\ge u_1(x)$ and $w(x)\ge u_2(x)$. Furthermore, for fixed $x_0\in G$ it holds $w(x_0)= u_1(x_0)$ or $w(x_0)= u_2(x_0)$, without loss of generality assume $w(x_0)= u_1(x_0)$. Then, from GCP we have $f(x_0)\le I(u_1,x_0)\le I(w,x_0)$, which concludes the proof as $x_0\in G$ is arbitrary.
    	
    \end{proof}

We conclude this discussion with the observation that Assumptions \ref{assume:Uniqueness} and  \ref{assume:Existence} together in fact imply that $I$ must actually be continuous.  We have opted to state Assumptions \ref{assume:Uniqueness} and \ref{assume:Existence} separately to emphasize that the comparison result Theorem~\ref{thmIntro:Comparison} is independent of the existence result Theorem~\ref{thmIntro:Existence} and requires only Assumption \ref{assume:Uniqueness} (however, note that Theorem~\ref{thmIntro:Existence} does indeed require both Assumptions  \ref{assume:Uniqueness} and \ref{assume:Existence}). However, the GCP in conjunction with other assumptions yields stronger conditions. For example, one can see the inequality in Assumption \ref{assume:Existence}~\eqref{assume:LSCBump} automatically holds for $x\neq x_0$ if $I$ also satisfies the GCP. We record here another consequence of combining the assumptions.

\begin{lemma}\label{lemIntro:IContinuous}
If $I$ satisfies Assumptions \ref{assume:Uniqueness}~\eqref{assume:GCP}-\eqref{assume:SubtractConstant} and \ref{assume:Existence}~\eqref{assume:UpperSemicon}-\eqref{assume:LSCBump}, then $I$ is continuous.

\end{lemma}

\begin{proof}[Proof of Lemma \ref{lemIntro:IContinuous}]
	
	We will only prove the continuity of $I$ at $u\equiv 0$, and the argument for a generic $u$ follows analogously.  As $I$ is u.s.c from Assumption~\ref{assume:Existence}~\eqref{assume:UpperSemicon}, we will only need to establish that $I$ is l.s.c.  To this end, assume that $\{u_k\}_{k=1}^\infty\subset C(G)$ with
	\begin{align*}
		u_k\to 0\ \text{in}\ C(G).
	\end{align*}
	Let $y\in G$ be fixed.  We will show that
	\begin{align*}
		\liminf_{k\to\infty} I(u_k,y)\geq I(0, y).
	\end{align*}
	To this end, let us define the constant, $c_k$, as
	\begin{align*}
	c_k\coloneqq \min\bigg(0,\min_{1\le j\le N}u_k(x_j)\bigg).
	\end{align*}
	We note that
	\begin{align*}
		c_k\leq 0\ \ \ \text{and}\ \ \ 
		\forall\ z\in G,\ c_k\leq u_k(z).
	\end{align*}
	Thus, we see that for the bump function, $b_y$, in Assumption \ref{assume:Existence}~\eqref{assume:LSCBump},
	\begin{align*}
		u_k - c_k \geq 0=(u_k-c_k)b_y \ \ \text{on }\ G\setminus \{y\}\ \text{and}\ \ 
		(u_k-c_k)(y) = (u_k(y)-c_k)b_y(y),
	\end{align*}
	which means the GCP is applicable in this case.  Furthermore, as $u_k\leq 0$, the inequality from Assumption \ref{assume:Uniqueness}~\eqref{assume:SubtractConstant} is reversed, hence
	\begin{align*}
		I(u_k,y)\geq I(u_k-c_k, y).
	\end{align*}
	Combining these and invoking the GCP, we see that
	\begin{align*}
		I(u_k,y)\geq I(u_k-c_k, y)\geq I((u_k-c_k)b_y, y).
	\end{align*}
	Thus, 
	\begin{align*}
		I(u_k,y)\geq \inf_{z\in G} \left\{ I((u_k-c_k)b_y, z)  - I(0, z) \right\}   +I(0, y).
	\end{align*}
	Invoking Assumption \ref{assume:Existence}~\eqref{assume:LSCBump}, we see that
	\begin{align*}
		\liminf_{k\to\infty} I(u_k,y) \geq I(0, y).
	\end{align*}

\end{proof}


\section{Proofs of Theorems \ref{thmIntro:Comparison} and \ref{thmIntro:Existence}}

	In this section, we present the proofs of Theorems \ref{thmIntro:Comparison} and \ref{thmIntro:Existence}.

	\begin{proof}[Proof of Theorem \ref{thmIntro:Comparison}]

		Let us start by considering the simpler case in which the inequality in the hypotheses of the theorem is strict.   That is to say, we assume that $u$, $v\in C(G)$ and that
		\begin{align*}
			I(u,x)> I(v,x)\ \ \text{for}\ \ x\in G\setminus\Gam.
		\end{align*}
		Let us argue by contradiction. Assume that our thesis is not true. This means that there exists $\tilde{x} \in G\setminus \Gamma$ such that
	    \begin{align*}
	    	(u-v)(\tilde{x})=\max_G(u-v)>\max_{\Gamma}(u-v),
	    \end{align*}
		and $(u-v)(\tilde{x})>0$ (otherwise, the conclusion is trivial).

		Let us define now the function
		\begin{align*}
			\tilde{u}\coloneqq  u-\tilde{c},\quad \tilde{c}\coloneqq (u-v)(\tilde{x})\geq 0.
		\end{align*}
		First, by Assumption \ref{assume:Uniqueness}~\eqref{assume:SubtractConstant} and the fact that $u$ is a strict subsolution, we get 
		\begin{align*}
			I(\tilde{u},\cdot)\ge I(u,\cdot)>I(v,\cdot)\ \ \text{on}\ \ G\setminus\Gam.
		\end{align*}
		Moreover, from the definition of $\tilde{u}$, we have $\tilde{u}(\tilde{x})=v(\tilde{x})$ and $\tilde{u}\le v$ on $G$. Thus, we can apply GCP to the functions $\tilde{u}$ and $v$, along with the previous inequality to obtain
		\begin{align*}
			I(v,\tilde x) < I(\tilde{u},\tilde{x})\le I(v,\tilde{x}),
		\end{align*}
		which yields a contradiction.

		Let us deal with the general case. We assume $u$ satisfies $I(u,\cdot)\geq I(v,\cdot)$ on $G\setminus\Gamma$ and the inequality is not necessarily strict.  We will assume that additionally, condition \eqref{assume:StrictSubSolPerturbation} in Theorem \ref{thmIntro:Comparison} is satisfied.  If it so happened that instead, condition \eqref{assume:StrictSuperSolPerturbation} was satisfied, then the argument is analogous.  We use the functions $u_k\in C(G)$ so that $I(u_k,\cdot)>I(u,\cdot)$ on $G\setminus\Gamma$. For any fixed $k$, the above proof implies that 
		\begin{align*}
		\max_G(u_k-v)_+\leq \max_\Gamma(u_k-v)_+,
		\end{align*}
		then taking $k\to \infty$ finishes the proof.

	\end{proof}

	\begin{proof}[Proof of Theorem \ref{thmIntro:Existence}]
		
		The proof follows the lines of the classical Perron's method in the continuous case. Let us define the set of admissible subsolutions, $\S$, as

		\begin{align*}
			\S\coloneqq 
			\{\psi:\,\,I(\psi,\cdot)\ge f\mbox{ on }G\setminus \Gamma, 
			\,\,\psi\le g \mbox{ on }\Gamma, 
			\text{ and }  \psi\geq U_{f,g,\sub}\ \text{on}\ G  \},
		\end{align*}
		 where $U_{f,g,\sub}$ is the subsolution that appears in the assumption of the theorem.  We define $u$ as the supremum
		\begin{align*}
			u(x)\coloneqq  \sup_{\psi\in \S}\psi(x),\quad x \in G,
		\end{align*}
		and we want to show that $u$ is the unique solution to \eqref{eqIntro:HJBGeneric}.   We note that once we have shown that $u$ is indeed a solution of \ref{eqIntro:HJBGeneric}, then the uniqueness is immediate from the comparison condition \eqref{eqIntro:ExistenceTheorem-ComparisonAssumption}, thus we need only verify that $u$ indeed solves \eqref{eqIntro:HJBGeneric}.

		We first focus on the subsolution property of $u$.  Fix $x_0 \in G\setminus\Gamma$, then we will show that $I(u,x_0)\geq f(x_0)$.  Let us take a sequence $\psi_n\in \S$ such that $\psi_n(x_0)\nearrow u(x_0)$. We first note that $\psi_n\le u$ for every $n$, by definition of $u$. Furthermore, since $\psi_n\in \S$, we know that $U_{f,g,\sub}\leq \psi_n$, while since $\psi_n\leq g$ on $\Gam$, by the comparison of sub and supersolutions, condition \eqref{eqIntro:ExistenceTheorem-ComparisonAssumption} in the Theorem, we see that  $\psi_n\le V_{f,g,\super}$, with $V_{f,g,\super}$ the supersolution that assumed to exist in the statement of the theorem. Recalling that $C(G)\cong\R^N$ this implies we can identify $\{\psi_n\}_{n=1}^\infty$ with a bounded sequence in $\R^N$, hence passing to a (not relabeled) subsequence we may assume there is some $\bar{\psi}\in C(G)$ such that $\psi_n\to\bar{\psi}$ pointwise on $G$. Note that $\bar{\psi}(x_0)=u(x_0)$ because $\psi_n(x_0)\nearrow u(x_0)$. Applying GCP, we then obtain $I(\bar{\psi},x_0)\le I(u,x_0)$.  Thus, combining this with the u.s.c. property of $I$ from Assumption~\ref{assume:Existence}~\eqref{assume:UpperSemicon}, we have

		\begin{align*}
			I(u,x_0)\geq
			I(\bar{\psi},x_0)\ge \limsup\limits_{n\to \infty}I(\psi_{n},x_0)\ge f(x_0),
		\end{align*}
		  which shows $u$ is a subsolution.

		Next, we want to demonstrate $u$ is a supersolution with the correct boundary values, i.e.
		\begin{align}\label{eqProofs:UIsSuperSolution}
			\begin{cases}
				I(u,\cdot)\le f& \mbox{on }G\setminus \Gamma,\\
				u= g & \mbox{on } \Gamma.
			\end{cases}
		\end{align}
		Let us argue the inequality in $G\setminus \Gam$ by contradiction, so we assume there is $x_0 \in G\setminus\Gamma$ and $I(u,x_0)> f(x_0)$. We take the bump function $b_{x_0}$ defined in Assumption \ref{assume:Existence}~\eqref{assume:LSCBump}, and we consider $u+tb_{x_0}$ with $t>0$.  Since 
		\begin{align*}
			u(x_0)+tb_{x_0}(x_0) > u(x_0),
		\end{align*}
		if we show that $u+tb_{x_0}\in \S$ then we will have a contradiction with the definition of $u$ as a pointwise supremum over $S$. As $x_0\in G\setminus\Gamma$, it is immediate that $u+tb_{x_0}\le g$ on $\Gamma$.

		It remains to prove $I(u+tb_{x_0},\cdot)\ge f$ on $G\setminus\Gamma$.  First, we will focus on $x=x_0$, and then we will address $x\not=x_0$ second. Since $I(u,x_0)>f(x_0)$, there exists a $\delta>0$ such that $I(u,x_0)>f(x_0)+\delta$. Then by Assumption \ref{assume:Existence}~\eqref{assume:LSCBump}, for all $t$ sufficiently small it holds 
		\begin{align*}
			\inf_{x \in  G\setminus\Gamma}(I(u+tb_{x_0},x)-I(u,x))\ge -\delta.
		\end{align*}
		Thus,
		\begin{align*}
			I(u+tb_{x_0},x_0)
			=I(u+tb_{x_0},x_0)- I(u,x_0) +I(u, x_0)\geq -\del+I(u, x_0)\geq f(x_0),
		\end{align*}
		and so the subsolution inequality is satisfied at $x=x_0$ when $t$ is small enough.  Now, assume that $x\not=x_0$.
We note that  $u+tb_{x_0}\geq u$ on all of $G$, and when $x\not=x_0$, we have $u(x)+tb_{x_0}(x)=u(x)$, thus by the GCP (Assumption \ref{assume:Uniqueness}~\eqref{assume:GCP}) applied at $x$, and the fact that we know $u$ is a subsolution, we obtain $I(u+tb_{x_0},x)\geq I(u,x)\geq f(x)$. In particular, this verifies that for any $t>0$ small enough, we have $u+tb_{x_0}\in \S$, yielding our desired contradiction.  Thus, we conclude that $I(u,\cdot)\leq f$ in $G\setminus\Gam$, i.e. that $u$ is a supersolution.  This yields that $I(u, \cdot)=f$ on $G\setminus \Gamma$.

Next, we verify the boundary data $u=g$ on $\Gam$.  By definition, as $u$ is a supremum over $\S$, and each $\psi\in \S$ satisfies $\psi\leq V_{f,g,\super}\leq g$ on $\Gam$, we naturally have  $u\leq g$ on $\Gam$. To see the opposite inequality, fix $x_\Gamma\in \Gamma$ and let $\psi\in \S$, and define
\begin{align*}
	\tilde\psi(x)
	=
	\begin{cases}
		\psi(x),&x\neq x_\Gamma,\\
		g(x_\Gamma),&x=x_\Gamma.
	\end{cases}
\end{align*}
we see that $\tilde\psi\geq \psi$ on $G$ and $\tilde\psi(x)=\psi(x)$ for any $x\in G\setminus \Gamma$, thus by the GCP  (Assumption \ref{assume:Uniqueness}~\eqref{assume:GCP}) we obtain 
\begin{align*}
	I(\tilde\psi, x)\geq I(\psi, x)\geq f(x)\ \  \text{for}\ \ x\in G\setminus \Gamma.
\end{align*}
Since $\tilde \psi\geq \psi\geq U_{f,g,\sub}$ on $G$ and $\tilde\psi\leq g$ on $\Gamma$, we see $\tilde\psi\in \S$, thus $u(x_\Gamma)\geq \tilde\psi(x_\Gamma)=g(x_\Gamma)$. 
		
	\end{proof}


\section{Examples}\label{sec:Examples}

	In this section, we provide some examples for which Assumptions \ref{assume:Uniqueness}, \ref{assume:Existence}, and the additional requirements of Theorems \ref{thmIntro:Comparison} and \ref{thmIntro:Existence} are satisfied.  One set of examples comes from \cite{Cald_Ette}, and we also show that those examples also have the min-max structure given in \eqref{eqIntro:MinMax} (in fact, only a max, as those examples are convex).  Finally, at the end of the section, we will show that for a large class of operators, the structure in \eqref{eqIntro:MinMax} is generic for our operators here as well as those in \cite{Cald_Ette}.


\subsection{Graph Laplacian for a Markov Chain}\label{subsec:Markov}

	For \eqref{eqIntro:HJBGeneric}, when $I$ happens to be linear and of a particular form, below, there is an explicit formula for the function, $u$, as an expected exit value of a corresponding Markov chain on $G$, when it hits $\Gam$ (see \eqref{eqMarkov:SolutionHJBProbGraphLaplace}).   We note that this is all standard, and we do not present anything new.  Rather, we point out this connection for its useful consequences (in later examples).  Specifically, if one wishes to uniquely solve \eqref{eqIntro:HJBGeneric} in the context of the generator of a Markov chain, then the existence of a solution comes from an explicit formula, and the uniqueness of solutions (which results from Assumption \ref{assume:Uniqueness}) is what needs to be checked.  Nonetheless, we will demonstrate that \emph{under certain restrictions on the corresponding Markov chain}, this class of linear operators satisfies both Assumptions \ref{assume:Uniqueness} and \ref{assume:Existence}.  Subsequently, we can use this construction to create barrier functions for other examples of equations.

	To this end, we assume that for $t\in\Natural\union\{0\}$, $X_t$ is a Markov chain on $G$ with transition density function, $K(x,y)$, on some probability space, $\Om$ with probability $\Prob$ and expectation, $\Expectation$.  That is to say, given $\Om$ with probability $\Prob$ and expectation $\Expectation$, we assume there exists a function
	\begin{align}\label{eqMarkov:DefOfK}
		K:G\times G\to [0, \infty),\ \ \text{with}\ \ \ \forall\ x\in G,\ \ 
		\sum_{y\in G} K(x,y)=1,
	\end{align}
	then it is standard that there exists a Markov chain $X_t$ satisfying
	\begin{align}\label{eqMarkov:KernelProbGraphLaplace}
		\Prob(X_{t+1}=y\mid X_t=x) = K(x,y).
	\end{align}
	The corresponding graph Laplacian for $X_t$ is defined using $K$, as
	\begin{align}\label{eqMarkov:ProbGraphLaplace}
		L(u,x)\coloneqq  \sum_{y\in G} (u(y)-u(x))K(x,y),
	\end{align}
	which is the generator of $X_t$. That is to say, conditioned upon $\Prob(X_{t_0}=x)=1$,
	\begin{align}\label{eqMarkov:Generator}
		&\mathbb{E}_x\left(u(X_{t_0+1})-u(X_{t_0})\right)=\sum_{y\in G}u(y)\mathbb{P}(X_{t_0+1}=y\mid X_{t_0}=x)-u(x)\nonumber\\
		&=\sum_{y\in G}u(y)K(x,y)-\sum_{y\in G}u(x)K(x,y)=L(u,x).
	\end{align}
	One of the many reasons the generator of $X_t$ is so important is that it allows for the following representation, which can be thought of the most basic instance of Dynkin's formula.   We do not include the definition of stopping time here, as this section is just meant to be a brief description of a useful example for elliptic equations on graphs.

	\begin{lemma}[Dynkin]\label{lemMarkov:Dynkin}
		
		If $w\in C(G)$, $X_t$ is as above, $L$ is its corresponding generator \eqref{eqMarkov:ProbGraphLaplace}, and $\tau$ is a stopping time that satisfies $\Expectation_x\tau <\infty$, then
		\begin{align}\label{eqMarkov:DynkinFormula}
			\Expectation_x(w(X_\tau)) = w(x) 
			+ \Expectation_x\left(   \sum_{t=0}^{\tau-1} L(w,X_t)   \right),
		\end{align}
	\end{lemma}

	\begin{proof}[Proof of Lemma \ref{lemMarkov:Dynkin}]
		The proof is a simple implementation of the argument that was used to obtain the generator in the first place.  Note the fact that $\Expectation_x\tau<\infty$ implies the quantities  in the calculation below are all well-defined:
			\begin{align*}
				\Expectation_x(w(X_\tau))
				&= w(x) + \Expectation_x\left(   \sum_{t=0}^{\tau-1} (w(X_{t+1})-w(X_t))   \right)\\
				&= w(x) + \Expectation_x\left(  \sum_{t=0}^{\tau-1}  
				\sum_{y\in G}\Prob(X_{t+1}=y \mid X_t=X_t)(w(y)-w(X_t))  \right)\\
				&= w(x) + \Expectation_x\left(   \sum_{t=0}^{\tau-1} \sum_{y\in G} K(X_t,y)(w(y)-w(X_t))   \right)\\
				&= w(x) + \Expectation_x\left(   \sum_{t=0}^{\tau-1} L(w,X_t)   \right).
			\end{align*} 
		
	\end{proof}

	Now, let $x$ be fixed and write $X_{x,t}$ to denote the process $X_t$ starting at $x$. We also denote by
	\begin{align}\label{eqMarkov:HittingTime}
		\tau_{x,\Gamma}\coloneqq  \min\{s\in\Natural\union\{0\}\mid  X_{x,s}\in \Gamma\},    	
	\end{align}
	the first hitting time of $X_{x,t}$ with respect to $\Gamma$. We remark that this is indeed a stopping time, and in the discrete setting, this time coincides with the exit time from $G\setminus \Gam$. Since $\Gamma$ is fixed throughout the work, we drop the explicit dependence of $\tau_{x,\Gamma}$ on $\Gamma$ and simply write $\tau_x$ instead.  \emph{In this setting, it is not always true that $K$, $G$, and $\Gam$ are such that $\tau_{x,\Gam}$ will be finite, and even have finite expectation.}  This is an assumption that needs to be made on $X_t$ in order for this example to be well defined.

	\begin{assumption}[Finite exit time]\label{assumeMarkov:FinteExitTimeExpectation}
		$K$, $G$, $\Gam$ are such that for each $x\in G\setminus\Gam$,
		
		\begin{align*}
			\Expectation_x\left( \tau_{x,\Gam} \right)<\infty.
		\end{align*}
	\end{assumption}

	\noindent
	Finally, we define the function $u$ as the expected running cost plus exit value of $X_t$ from $G\setminus \Gam$:
	\begin{align}\label{eqMarkov:SolutionHJBProbGraphLaplace}
		u(x)\coloneqq  \mathbb{E}_{x}\bigg(\sum_{t=0}^{\tau_x-1}f(X_{x,t})+g(X_{x,\tau_x})\bigg).
	\end{align}
	We note that thanks to Assumption \ref{assumeMarkov:FinteExitTimeExpectation}, $u$ is well defined on $G$. Additionally, if $X_0\in\Gam$, there should be no contribution to the running cost, thus we adopt the convention that 
	\begin{align}\label{eqMarkov:SumConventionStartInGam}
		\text{for}\ \ i>j,\ \ 
		\sum_{t=i}^{j}f(X_{x,t})=0.
	\end{align}

	The reason we have chosen to address this example is that $u$ solves \eqref{eqIntro:HJBGeneric}, for $-f$, when $I(u,x) =L(u,x)$, as in \eqref{eqMarkov:ProbGraphLaplace}.

	\begin{proposition}\label{propMarkov:ExistenceProbGraphLaplace}

		If Assumption \ref{assumeMarkov:FinteExitTimeExpectation} is satisfied, $u$ is defined as in \eqref{eqMarkov:SolutionHJBProbGraphLaplace}, and $L$ is defined in \eqref{eqMarkov:ProbGraphLaplace}, then $u$ is a solution to
		\begin{align}\label{eqMarkov:LinearEqProbGraphLaplace}
			\begin{cases}
				L(u,\cdot)=-f& \mbox{in }G\setminus \Gam,\\
				u=g & \mbox{on } \Gamma.
			\end{cases}
		\end{align}
		Furthermore, if Assumption \ref{assumeMarkov:FinteExitTimeExpectation} is satisfied, $w\in C(G)$ exists and satisfies \eqref{eqMarkov:LinearEqProbGraphLaplace}, then $w=u$ where $u$ is as in \eqref{eqMarkov:SolutionHJBProbGraphLaplace}.

	\end{proposition}

	\begin{proof}[Proof of Proposition \ref{propMarkov:ExistenceProbGraphLaplace}]
		We start with the first implication.  We note that by the sum convention given in \eqref{eqMarkov:SumConventionStartInGam}, if $x\in\Gam$, then $u(x)=g(x)$.

		It remains to show $L(u,x)=-f(x)$ if $x\in G\setminus \Gam$.  Since $x\not\in\Gam$, we know that $\tau_x\geq 1$. Here we use that fact that we can start $X_0=x$, run for one step, and then recalculate the remaining probabilities.  By the Markov property, we can restart the process at its one-step location.  If $X_{x,1}=z$, let us abuse notation, and call the restarted process $X_{z,t}$ with $X_{z,0}=z$, and its exit time will be denoted as $\tau_z$.   That is, we have 
		\begin{align*}
			&\Expectation_x\left(
			\sum_{t=0}^{\tau_x-1}f(X_{x,t})+g(X_{x,\tau_x})
			\right)
			= \Expectation_x\left( f(X_{x,0}) + 
			\sum_{t=1}^{\tau_x-1}f(X_{x,t})+g(X_{x,\tau_x})
			\right)  \\
			&= \sum_{z\in G}
			\Prob(X_{x,1}=z\mid X_{x,0}=x) 
			\left( f(x) +
			 \sum_{t=0}^{\tau_z-1}f(X_{z,t})+g(X_{z,\tau_z})
			\right),
		\end{align*}
		which specifically shows that
		\begin{align*}
			u(x) = f(x)\sum_{z\in G} K(x,z) + \sum_{z\in G} K(x,z) u(z)
			= f(x) + \sum_{z\in G} K(x,z) u(z).
		\end{align*}
		Rearranging, and recalling the defintion of $L(u,x)$ from \eqref{eqMarkov:ProbGraphLaplace}, we obtain
		\begin{align*}
		-f(x)=-\sum_{z\in G}K(x,z)u(x)+\sum_{z\in G}K(x,z)u(z)=L(u,x).
		\end{align*}

Next, we will address the second implication.  To this end, we assume we have a $w\in C(G)$ that solves \eqref{eqMarkov:LinearEqProbGraphLaplace}.  Let $X_t$ be the Markov chain corresponding to $L$.   We note that by assumption, $\Expectation_x\tau_x<\infty$.
 Plugging $\tau_x$ into Dynkin's formula, rearranging and using the equation, we see that for the process, $X_t$, above,
\begin{align*}
	w(x) = \Expectation_x\left(  -\sum_{t=0}^{\tau_x-1} L(w,X_t)  + w(X_{\tau_x})   \right)
	= \Expectation_x\left(   \sum_{t=0}^{\tau_x-1}f(X_t)  + g(X_{\tau_x})   \right) 
    =u(x).
\end{align*}	
		
	\end{proof}

	We can see that Proposition \ref{propMarkov:ExistenceProbGraphLaplace} gives the existence of solutions of \eqref{eqIntro:HJBGeneric}, and so Assumption \ref{assume:Existence} is technically irrelevant here.  Nonetheless, we will show that for $I=L$ in \eqref{eqMarkov:ProbGraphLaplace}, all of the relevant assumptions are indeed satisfied.
	
	\begin{proposition}[Assumptions for uniqueness]\label{propMarkov:L-SatisfiesAssumptionsUniqueness}
		
		$L$ satisfies Assumption \ref{assume:Uniqueness} and Theorem \ref{thmIntro:Comparison}~\eqref{assume:StrictSubSolPerturbation}.
	\end{proposition}
	
	\begin{proof}[Proof of Proposition \ref{propMarkov:L-SatisfiesAssumptionsUniqueness}]
		
		The non-negativity of $K$ and the structure of $L$ in \eqref{eqMarkov:ProbGraphLaplace} immediately give the GCP and the behavior with respect to subtracting positive constants (Assumptions \ref{assume:Uniqueness}~\eqref{assume:GCP} and \eqref{assume:SubtractConstant}).  The strict subsolution perturbation (i.e. condition Theorem \ref{thmIntro:Comparison}~\eqref{assume:StrictSubSolPerturbation}) follows from Lemma \ref{lemIntro:StrictSubsolPerturbation}~\eqref{assumptionLemIntro:BarrierFunction} and the linearity of $L$, by defining $\phi_u$ by \eqref{eqMarkov:SolutionHJBProbGraphLaplace} with $f=-1$ and applying Proposition \ref{propMarkov:ExistenceProbGraphLaplace}.  We note that in this case, the same barrier in Lemma \ref{lemIntro:StrictSubsolPerturbation}~\eqref{assumptionLemIntro:BarrierFunction} works for all choices of possible subsolutions.
		
	\end{proof}

	\begin{proposition}[Assumptions for existence]\label{propMarkov:L-SatisfiesAssumptionExistence}
		
		$L$ satisfies Assumption \ref{assume:Existence}, and assumptions \eqref{eqn: special sub/super} and  \eqref{eqIntro:ExistenceTheorem-ComparisonAssumption} in Theorem \ref{thmIntro:Existence}.
	\end{proposition}
	
	\begin{proof}[Proof of Proposition \ref{propMarkov:L-SatisfiesAssumptionExistence}]
		
		First, we note that by construction, $L$ is a linear and bounded operator on $C(G)$, and thus continuous, in particular the u.s.c. Assumption \ref{assume:Existence}~\eqref{assume:UpperSemicon} holds.  The continuity of $L$ also establishes the inequality required for the bump function in Assumption \ref{assume:Existence}~\eqref{assume:LSCBump}, however, it can be useful to see the inequality achieved by a different argument, which we mention here.  By the linearity and explicit structure of $L$, we see that for $x_0$ fixed,
		\begin{align*}
			L(u+tb_{x_0},x) - L(u,x)
			= tL(b_{x_0},x)
			= 
			\begin{cases}
				 t\sum_{y\in G\setminus \{x_0\}} (-1)K(x_0,y) \geq -t, &\text{if}\ x=x_0,\\
				 tK(x,x_0), &\text{if}\ x\not=x_0.
			\end{cases}
		\end{align*}
		We note the first inequality is from that fact that
		\begin{align*}
			\sum_{y\in G}K(x_0,y) = 1,\ \ \text{and so}\ \ 
			0\leq \sum_{y\in G\setminus \{x_0\}} K(x_0,y) \leq 1.
		\end{align*}
		Hence, for $t>0$,
		\begin{align*}
			L(u+tb_{x_0},x) - L(u,x)\geq -t,
		\end{align*}
		and so Assumption \ref{assume:Existence}~\eqref{assume:LSCBump} holds.  For the specific sub/supersolutions required in Theorem \ref{thmIntro:Existence}, if $f$ and $g$ are given as in \eqref{eqIntro:HJBGeneric}, then for $U_{f,g,\sub}$ and $V_{f,g,\super}$ we can take the functions defined in \eqref{eqMarkov:SolutionHJBProbGraphLaplace} with $g=g$ and $f$ replaced by $\max_G(f)$ and $\min_G(f)$ respectively; Proposition \ref{propMarkov:ExistenceProbGraphLaplace} yields the desired sub/supersolution properties. Finally, by Proposition \ref{propMarkov:L-SatisfiesAssumptionsUniqueness} above, we may apply the comparison Theorem \ref{thmIntro:Comparison} to see that \eqref{eqIntro:ExistenceTheorem-ComparisonAssumption} is satisfied.
	\end{proof}

	\begin{rem}
		We note that the structure of $L$ in \eqref{eqMarkov:ProbGraphLaplace}, neglecting the unit mass of $K(x,\cdot)$, as well as those operators as in \eqref{eqIntro:IGraphLaplace}, is generic for the linear case.  In fact, any \emph{linear} operator on $C(G)$, say $I$, that also has the GCP must necessarily have corresponding functions, $c: G\to\real$ and $a:G\times G\to\real$, with $a\geq 0$, and
		\begin{align*}
			L(u,x) = c(x)u(x) + \sum_{y\not=x} \left(u(y)-u(x)\right) a(x,y).
		\end{align*}
		(This basic observation was demonstrated and used for the min-max structure of elliptic operators in \cite[Lemma 2.1]{GuSc-2019MinMaxNonlocalCALCVARPDE}.)
	\end{rem}

	
	\subsection{A Control problem for Markov chains (the Bellman equation)}\label{subsec:ControlProblem}

	A class of nonlinear and convex examples of $I$ can be obtained by considering a control problem related to the Markov chain in the previous subsection.  We note there are many types of control problems for Markov chains.  The control problem we use allows a participant to choose from a set of various transition densities (like those $K$ in \eqref{eqMarkov:DefOfK}) to drive the process to a desired outcome of running cost and boundary value.  To this end, we introduce some further notation to set up the control problem in question.   What appears here is an example of a Bellman equation from optimal control.  Again, as in the previous subsection, the framework presented here is not new, but is included as an example of an application of our results.
	
	Let  $\A$ be  a collection of indices which we use to label the set of admissible transition densities associated to the control of the Markov chain: 
	\begin{align*}
		\forall\ i\in\A,\ K^i\ \text{satisfies}\ (\ref{eqMarkov:DefOfK}).
	\end{align*}
	For this type of problem, a \emph{control} is, at each element $x\in G$, a choice of $i\in\A$ which gives the next step of the chain dictated by $K^i$.  In other words, to try to emphasize this structure, equivalently a control is a choice of function, $\al:G\to\A$, which, at each $x$, selects the transition density $K^{\alpha(x)}$.  Then for this fixed choice of $\al$, a Markov process $X^\al$ can be constructed satisfying the transition law
	\begin{align*}
		\Prob(X^\al_{t+1}=y \mid X^\al_{t}=x) = K^{\al(x)}(x,y).
	\end{align*}
	 We note that by this choice of construction, all of the resulting Markov chains, $X^\al_t$, are all defined on the same probability space, $\Om$, with the same background probability measure, $\Prob$, and same expectation, $\Expectation$.  For ease of notation, we will write 
	 \begin{align*}
	K^{\al}(x,y) \coloneqq  K^{\al(x)}(x,y);
\end{align*}
	 we then see that the assumption that at each $x$, $K^{\al(x)}$, satisfies (\ref{eqMarkov:DefOfK}), this transition kernel, $K^\al$, indeed defines a Markov chain.

To set some notation, we already have the notion of controls (functions from $G$ to $\al$), and we need to collect the set of all possible kernels that can result from this construction.  We take the following:
\begin{align}\label{eqControl:ControlSetAndKernelSet}
	\C \coloneqq  \{\al\ |\ \al:G\to\A\}
	\ \ \ \text{and}\ \ \ 
	\K \coloneqq  \{ K^\al\ |\ \al\in\C \}.
\end{align}
Note the set $\mathcal{K}$ is in general much larger than the index set $\mathcal{A}$: for example, if there are two kernels $K^1$ and $K^2$, with $\A = \{1,2\}$, because we consider all possible controls the class $\mathcal{K}$ contains $\lvert G\rvert^2$ kernels.

	Now, with this notation, the associated running cost, respectively terminal cost, are the function $f$, respectively $g$, from the previous section in \eqref{eqMarkov:SolutionHJBProbGraphLaplace}. As with the Markov chain of the previous example, we will need an assumption on the class of admissible transition densities that ensure the controlled process always exits the domain in finite time.  Given $\al\in\C$, we denote the process starting at $x$ for $t=0$ by $X^\al_{x,t}$, and the exit time as,
	\begin{align*} 
		\tau^\al_{x,\Gamma}\coloneqq  \min\{s\in\Natural\union\{0\}\mid  X^\al_{x,s}\in \Gamma\}.
	\end{align*}  
	Given a specific control, $\al$, that would give an associated cost corresponding to $X^\al_{x,t}$ and its exit time $\tau^\al_x$, denoted by $u^\al$ as in \eqref{eqMarkov:SolutionHJBProbGraphLaplace}, which for subsequent ease of notation, we denote by
	\begin{align*}
		u^{\alpha}(x)\coloneqq \mathbb{E}_{x}\bigg(\sum_{t=0}^{\tau^{\alpha}_x-1}f(X^{\alpha}_{x,t})+g(X^{\alpha}_{x,\tau^{\alpha}_x})\bigg).
	\end{align*}  
	The goal is now to find a least cost among all possible choices of a control. 
	In order that we can define a minimal cost, we make the following assumption on $\tau^\al$.
	
\begin{assumption}\label{assumeControl:FiniteExitTime}
	
	The sets, $\C$ and $\K$ in (\ref{eqControl:ControlSetAndKernelSet}), are such that for each $x\in G\setminus\Gam$,
	\begin{align*}
		\sup_{\al\in\C}\Expectation_x(\tau^\al_{x,\Gam})<\infty.
	\end{align*}
	
\end{assumption}
This assumption says that the collection of all admissible kernels, $\K$, is ``non-degenerate'' in the sense that no matter which control is chosen, the resulting Markov chain will always exit with finite expected exit time.
	Under this assumption, we can define the minimum cost as
	\begin{align}\label{eqControl:SolutionControlProblem}
		u(x)\coloneqq  \inf_{\alpha\in \C}\mathbb{E}_{x}\bigg(\sum_{t=0}^{\tau^{\alpha}_x-1}f(X^{\alpha}_{x, t})+g(X^{\alpha}_{x, \tau^{\alpha}_x})\bigg)=\inf_{\alpha\in \C}u^{\alpha}(x).
	\end{align}
	We point out that each $u^{\alpha}$ is the unique solution of \eqref{eqMarkov:LinearEqProbGraphLaplace} with $K^\al$ associated to $X^{\alpha}$.  Now, we can demonstrate the equation that is satisfied by the infimum of these functions, which is $u$, above.  This is called the Bellman equation for the control problem above.

	\begin{theorem}[Bellman equation]\label{thmControl:BellmanEq}
		
		If $\K$ satisfies Assumption \ref{assumeControl:FiniteExitTime}, for all $\al\in\C$, $K^{\alpha}(x,x)=0$ for all $x\in G$, and $u$ is defined as in \eqref{eqControl:SolutionControlProblem}, then, $u$ solves
		\begin{align}\label{eqControl:BellmanEquation}
			\begin{cases}
				\displaystyle \inf_{\alpha\in \C}
				\left( \sum_{y\in G}(u(y)-u(x))K^{\alpha}(x,y) \right)
				= -f(x)\ &\text{in}\ G\setminus \Gam\\
				\displaystyle u=g\ &\text{on}\ \Gam
			\end{cases}
		\end{align}

	\end{theorem}

\begin{rem}
	We note that typically in the literature, the previous theorem is presented with a second part (similar to the second part in Proposition \ref{propMarkov:ExistenceProbGraphLaplace}), which is called verification.  It states that under certain assumptions of $u^\al$ being finite, and the existence of certain controls, that any solution of \eqref{eqControl:BellmanEquation} must have the form given in \eqref{eqControl:SolutionControlProblem}.  We do not need that second part for our investigations, and so we omit it for the sake of brevity.
\end{rem}

	Before giving the proof, we introduce some notation.  We let the operator $I$ be defined by
	\begin{align}\label{eqControl:BellmanOperatorInf}
		I(u,x)\coloneqq 
		\inf_{\alpha\in \C}\sum_{y\in G}(u(y)-u(x))K^{\al}(x,y),
	\end{align} 
	and given $K\in\K$, we have the linear operator as in subsection \ref{subsec:Markov}
	\begin{align}\label{eqControl:LinearOperatorLk}
		L_K(u,x) \coloneqq  \sum_{y\in G}(u(y)-u(x))K(x,y).
	\end{align} 
	For ease of notation, depending upon context, we can write $I$ in two different ways,
	\begin{align*}
		I(u,x) = \inf_{\alpha\in \C}\sum_{y\in G}(u(y)-u(x))K^{\alpha}(x,y)= \inf_{K\in\K} L_K(u,x),
	\end{align*}	
	where we recall that we are writing $K^\al(x,y)=K^{\al(x)}(x,y)$.

	\begin{proof}[Proof of Theorem \ref{thmControl:BellmanEq}]
		
		According to the same argument used to obtain the generator for $X_t$ in \eqref{eqMarkov:Generator}, we can rewrite
		\begin{align*}
				 u(x)&= \inf_{\alpha\in \C}\sum_{y\in G}K^{\alpha}(x,y)\mathbb{E}_{y}\bigg(\sum_{t=0}^{\tau^{\alpha}_{y}-1}f(X^{\alpha}_{y,t})+f(x)+g(X^{\alpha}_{y, \tau^\alpha_{y}})\bigg)\\
				&=\inf_{\alpha\in \C}\sum_{y\in G}K^{\alpha}(x,y)(f(x)+u^{\alpha}(y))\\
				&=f(x)+\inf_{\alpha\in \C}\sum_{y\in G}K^{\alpha}(x,y)u^{\alpha}(y).
		\end{align*}
		Moreover, by definition  of $u$ in \eqref{eqControl:SolutionControlProblem}, it holds $u^{\alpha}(y)\ge u(y)$ for any $\alpha\in \C$. Thus, we obtain
		\begin{align*}
			u(x)\ge f(x)+\inf_{\alpha\in \C}\sum_{y\in G}K^{\alpha}(x,y)u(y),
		\end{align*}
		and since $K^\alpha$ satisfies \eqref{eqMarkov:DefOfK} we obtain
		\begin{align*}
			-f(x)\ge \inf_{\alpha\in \C}\sum_{y\in G}K^{\alpha}(x,y)(u(y)-u(x)).
		\end{align*}

		To prove the opposite inequality, let $\ep>0$ be given, then note that for any $y\in G\setminus\Gam$, in which $y$ will be a new starting point for the Markov chain, by definition of (\ref{eqControl:SolutionControlProblem}), there exists $\alpha^*_y\in\C$ such that
		\begin{align*}
			u(y) + \ep \geq   \mathbb{E}_{y}\bigg(\sum_{t=0}^{\tau^{\alpha^*_y}_{y}-1}f(X^{\alpha^*_y}_{y,t})+g(X^{\alpha^*_y}_{y,\tau^{\alpha^*_y}_{y}})\bigg).
		\end{align*}
		Now fix $x\in G$, then there exists a control $\beta^*_x\in \C$, for which
		\begin{align*}
			\ep + \inf_{\alpha\in \C}\sum_{y\in G}K^{\alpha}(x,y)(u(y)-u(x))\geq \sum_{y\in G}K^{\beta^*_x}(x,y)(u(y)-u(x)).
		\end{align*}
		We define a new control $\tilde \al_x$ by
		\begin{align*}
			\tilde{\alpha}_x(y)\coloneqq 
			\begin{cases}
				\beta^*_x(x)&\mbox{if }y=x,\\
				\alpha^*_y(y)&\mbox{if }y\not=x,
			\end{cases}
		\end{align*}
		which we note satisfies $\tilde{\alpha}_x\in \C$.  By restarting the chain at $y$ after the first jump using the kernel for $\tilde \al_x$, we have $X^{\tilde \al_x}_{y,t}$ with $X^{\tilde\al_x}_{y,0}=y$, and recalling our assumption that $K^{\alpha(x)}(x,x)=0$ for any $x\in G$, we obtain
		\begin{align*}
			u(x)\le u^{\tilde{\alpha}_x}(x)&=\sum_{y\in G}K^{\tilde{\alpha}_x(x)}(x,y)\mathbb{E}_{y}\bigg(\sum_{t=0}^{\tau^{\tilde{\alpha}_x}_{y}-1}f(X^{\tilde{\alpha}_x}_{y,t})+f(x)+g(X^{\tilde{\alpha}_x}_{y,\tau^{\tilde\alpha_x}_{y}})\bigg)\\
			&=\sum_{y\in G}K^{\beta^*_x(x)}(x,y)\mathbb{E}_{y}\bigg(\sum_{t=0}^{\tau^{\alpha^*_y}_y-1}f(X^{\alpha^*_y}_{y,t})+f(x)+g(X^{\alpha^*_y}_{y,\tau^{\alpha^*_y}_y})\bigg)\\
			&\leq \sum_{y\in G}K^{\beta^*_x(x)}(x,y)(f(x)+u(y)+\ep),
		\end{align*}
		After rearranging, using the definition of $\beta^*_x$ and that $\sum_y K^{\beta^*_x}(x,y)=1$, we see that
		\begin{align*}
			-f(x)\le \inf_{\alpha\in \C}\sum_{y\in G}K^{\alpha(x)}(x,y)(u(y)-u(x)) + 2\ep.
	    \end{align*}
		Since $\ep>0$ was arbitrary, we conclude the desired inequality.

	\end{proof}

	Again, as in the previous subsection, at least under the special conditions Assumption~\ref{assumeControl:FiniteExitTime} for which $u$ is always finite, Theorem \ref{thmControl:BellmanEq} gives an explicit construction for the existence of solutions.  Below, we will show that $I$ defined in \eqref{eqControl:BellmanOperatorInf} satisfies Assumptions \ref{assume:Uniqueness} and \ref{assume:Existence}.  Before proving this, we note that associated to such $I$, there is a notion of minimal and maximal operators: 
	\begin{align*}
		\M_{\K}^-(u,x) \coloneqq  \inf_{K\in\K}\left( L_K(u,x) \right)=I(u, x)
		\ \ \ \text{and}\ \ \ 
		\M_\K^+(u,x) \coloneqq  \sup_{K\in\K}\left(  L_K(u,x)   \right),
	\end{align*}
	and it is straightforward to verify that
	\begin{align}\label{eqControl:MinimalMaximalInequalities}
		\forall\ u,v\in C(G),\ \ 
		\M_\K^-(u-v,x) 
		\leq I(u,x) - I(v,x)
		\leq \M_\K^+(u-v,x).
	\end{align}
	Furthermore, since we can identify $C(G)$ with $\real^N$ and we will always work in a situation for which $\M_\K^\pm$ are locally bounded, by their concavity, respectively convexity, $\M_\K^\pm$  are also continuous as operators from $C(G)$ to $C(G)$.  These extremal operators can be useful in verifying the assumptions related to uniqueness and existence.
	
	The operators $\M^\pm_\K$ play a fundamental role in well-posedness of the control problem.  One obstruction to well-posedness is that depending on the set of admissible transition densities, $\K$, there could exist a control such that the underlying Markov chain stays inside $G\setminus\Gam$ for arbitrarily long times.  In the event that $f<0$, this could yield that $u$ in \eqref{eqControl:SolutionControlProblem} is not finite (or is identically $-\infty$).  This phenomenon is avoided, for example, by Assumption \ref{assumeControl:FiniteExitTime}, above.  Another possible assumption to avoid this, is to \emph{assume} that $\K$ is such that there exists a $\phi\in C(G)$ satisfying
	\begin{align}\label{eqControl:MinimalEqSolution}
		\begin{cases}
			\M^-_{\K}(\phi, \cdot) = 1\ &\text{in}\ G\setminus\Gam\\
			\phi=0\ &\text{on}\ \Gam.
		\end{cases}
	\end{align}
	It turns out that these two assumptions are equivalent, which we note here.

	\begin{lemma}\label{lemControl:FiniteExitEquivSolBellmanRHS1}
		There exists $\phi$ solving \eqref{eqControl:MinimalEqSolution} if and only if Assumption \ref{assumeControl:FiniteExitTime} is satisfied.
	\end{lemma}

	\begin{proof}[Proof of Lemma \ref{lemControl:FiniteExitEquivSolBellmanRHS1}]
		The fact that Assumption \ref{assumeControl:FiniteExitTime} implies the existence of $\phi$ in \eqref{eqControl:MinimalEqSolution} is immediate by applying Theorem \ref{thmControl:BellmanEq} with $f=-1$.

		For the reverse direction, let $\phi$ be as in \eqref{eqControl:MinimalEqSolution}.  We wish to show that 
		\begin{align*}
			\sup_{x}\sup_{\al\in\C}\Expectation_x\left( \tau^\al_x \right) < \infty.
		\end{align*}
		To this end, let $\al\in\C$ be any control, $K^\al$ the corresponding transition kernel, and let $X_x^\al$ be the corresponding process started at $x$.  We note that by construction of $\M_\K^-$ and the equation for $\phi$, for $y\in G\setminus\Gam$
		\begin{align*}
			L_{K^\al}(\phi,y)\geq \M_\K^-(\phi,y) = 1.
		\end{align*}
		Since we are not assuming a priori that $\tau^\al_x$ is finite, for $k\in\Natural$ we define the new stopping time,
		\begin{align*}
			\tau_k = \min\{ k, \tau^\al_x\}.
		\end{align*}
		We may apply Dynkin's formula (Lemma \ref{lemMarkov:Dynkin}) to $\tau_k$ with $w=\phi$, to see that
		\begin{align*}
			\phi(x) &= \Expectation_x\left(   -\sum_{t=0}^{\tau_k-1} L_{K^\al}(\phi,X_t^\al) + \phi(X^\al_{\tau_k})   \right)\\
			&\leq \Expectation_x\left(  - \sum_{t=0}^{\tau_k-1} 1  + \phi(X^\al_{\tau_k}) \right)\\
			&\leq  \Expectation_x(-\tau_k)+\norm{\phi}_{C(G)}.
		\end{align*}
		Letting $k\to\infty$, we see that 
		\begin{align*}
			\phi(x) \leq \norm{\phi}_{C(G)} + \liminf_{k\to\infty} \Expectation_x (-\tau_k),
		\end{align*}
		thus
		\begin{align*}
			\limsup_{k\to \infty} \Expectation_x(\tau_k)
			\leq -\phi(x) + \norm{\phi}_{C(G)}
			\leq 2\norm{\phi}_{C(G)}.
		\end{align*}
		We note that $\tau_k\in\Natural\union\{0\}$, and for any $\omega\in \Omega$, the sequence $\{\tau_k(\omega)\}_{k=1}^\infty$ is pointwise increasing and
		\begin{align*}
			\lim_{k\to \infty} \tau_k(\om) = \tau^\al_x(\om)\in [0,+\infty].
		\end{align*}
		Thus, by the monotone convergence theorem we see that
		\begin{align*}
			\Expectation_x \left(  \tau^\al_x  \right)\leq 2\norm{\phi}_{C(G)}.
		\end{align*}
		As the right hand side is independent of $x$ and $\al$, we have attained the desired goal.

	\end{proof}

	Under the finite exit time Assumption \ref{assumeControl:FiniteExitTime}, $I$ will indeed satisfy the assumptions given in Section \ref{subsec:Assumptions}.

	\begin{proposition}[Assumptions for uniqueness]\label{propControl:SatisfiesAssumptionUniqueness}
		
		If Assumption \ref{assumeControl:FiniteExitTime} is satisfied, then
		$I$ defined by \eqref{eqControl:BellmanOperatorInf} satisfies Assumption \ref{assume:Uniqueness} and assumption Theorem \ref{thmIntro:Comparison}~\eqref{assume:StrictSubSolPerturbation}.
	\end{proposition}
	
	\begin{proof}[Proof of Proposition \ref{propControl:SatisfiesAssumptionUniqueness}]
		For the GCP (Assumption \ref{assume:Uniqueness}~\eqref{assume:GCP}), we note that it is immediate that for each $\al\in\C$, $L_{K^\al}$ has the GCP.  Since taking the infimum over $\al\in\C$ is a monotone operation, we see that also $I$ enjoys the GCP.  Since $I$ is an infimum over linear operators with the structure in \eqref{eqControl:LinearOperatorLk}, we obtain that for any constant function $c$ on $G$,
		\begin{align*}
			I(u-c,x) = I(u,x),
		\end{align*}
		and hence Assumption \ref{assume:Uniqueness}~\eqref{assume:SubtractConstant} is satisfied.  Finally we verify assumption in Theorem \ref{thmIntro:Comparison}~\eqref{assume:StrictSubSolPerturbation} of perturbing to strict subsolutions.  Indeed, by Lemma \ref{lemControl:FiniteExitEquivSolBellmanRHS1}, there exists $\phi$ that solves \eqref{eqControl:MinimalEqSolution}, then using \eqref{eqControl:MinimalMaximalInequalities}, and that $\M_\K^-$ is positively $1$-homogeneous, we see that
		\begin{align*}
			I(u+t\phi, \cdot) - I(u, \cdot) \geq \M_\K^-(t\phi, \cdot) = t\cdot 1>0.
		\end{align*}
		Again, we note, as in the linear case, this same choice of $\phi$ works for all possible $u$.
	\end{proof}

	\begin{proposition}[Assumptions for existence]\label{propControl:SatisfiesAssumptionExistence}
		
		$I$ as in \eqref{eqControl:BellmanOperatorInf} satisfies Assumption \ref{assume:Existence} and assumptions \eqref{eqn: special sub/super} and \eqref{eqIntro:ExistenceTheorem-ComparisonAssumption} in Theorem \ref{thmIntro:Existence}. 
	\end{proposition}
	
	\begin{proof}[Proof of Proposition \ref{propControl:SatisfiesAssumptionExistence}]
		We note that $I$ is concave, and by the fact that $C(G)$ is equivalent to $\real^N$, $I$ is a finite, concave function on $\real^N$, and hence continuous by \cite[Corollary 10.1.1]{RockafellarBook1970}.  In particular, the u.s.c. Assumption \ref{assume:Existence}~\eqref{assume:UpperSemicon} is satisfied.  Also, by continuity, Assumption \ref{assume:Existence}~\eqref{assume:LSCBump} is immediate, but, just as above for the linear case, we will note a different argument pertaining to the perturbation by the bump function.  Calculating as in the proof of Proposition \ref{propMarkov:L-SatisfiesAssumptionExistence}, we see that for each $K\in\K$,
		\begin{align*}
			L_K(tb_{x_0},x)\geq -t.
		\end{align*}
		Thus, using \eqref{eqControl:MinimalMaximalInequalities}, we see that
		\begin{align*}
			I(u+tb_{x_0},x) - I(u,x) 
			\geq \M_\K^-(tb_{x_0},x)\geq -t.
		\end{align*}
		As in the linear case, for the sub/supersolutions $U_{f,g,\sub}$ and $V_{f,g,\super}$, given $f$ and $g$ as in \eqref{eqIntro:HJBGeneric} we may apply Theorem \ref{thmControl:BellmanEq} with $g=g$ and $f$ in the theorem replaced by the constant functions $\max_G(f)$ and $\min_G(f)$ respectively. Also Proposition \ref{propControl:SatisfiesAssumptionUniqueness} allows us to apply the comparison Theorem \ref{thmIntro:Comparison} to see that \eqref{eqIntro:ExistenceTheorem-ComparisonAssumption} is satisfied.
		
	\end{proof}


\subsection{Graph eikonal and $p$-eikonal operators}\label{subsec:Eikonal}

		In this section, investigate two operators that appeared in \cite{BungertCalderRoith-2023UnifRatesLipLearningOnGraphs} and \cite{Cald_Ette}, and show that they fit within our framework here. We endeavor to keep notation similar to that in \cite{Cald_Ette}, but may deviate at certain points. In this section, we assume the vertices of $G$ are connected by edges with associated weights, and write
		\begin{align*}
			w(x,y)\ \text{is the edge weight from}\ x\ \text{to}\ y,
		\end{align*}
		for convenience, when it is clear from context, we will shorten this to
		\begin{align}\label{eqEikonal:EdgeWeightNotation}
			w_{x,y} \coloneqq  w(x,y)\ \ \ \text{or}\ \ \ 
			w_{ij} \coloneqq  w(x_i,x_j).
		\end{align}
		Following \cite{Cald_Ette} the graph eikonal operator is given as
		\begin{align*}
			H_e(u,x) \coloneqq  \max_{y\in G}\left( w_{y,x}(u(x)-u(y)) \right),
		\end{align*}
		and for $1\leq p<\infty$, the $p$-eikonal operator is
		\begin{align*}
			H_p(u, x)\coloneqq  \sum_{y\in G}\frac{1}{p}w_{y,x}((u(x)-u(y))_+)^p.
		\end{align*}
		In what follows, for simplified notation, we may sometimes write $u_i$ for $u(x_i)$.  Here we have added the factor of $\frac{1}{p}$ for some ease of analysis, this is an immaterial change as we are not concerned with studying limits as $p\to\infty$.

		\begin{rem}
			As in Section \ref{sec:Intro}, recall the sign convention in the GCP of Definition \ref{defIntro:GCP}, as compared to the convention in \cite[Definition 2]{Cald_Ette} and \cite[Proposition 3]{Cald_Ette}.  Indeed, if $u-v$ attains a zero maximum at $x_0$, then our requirement on $I$ says $I(u,x_0)\leq I(v,x_0)$, whereas if $H$ is an example from \cite[Section 2]{Cald_Ette}, it would satisfy $H(u,x_0)\geq H(v,x_0)$. 
		\end{rem}

		In light of this observation, we see the operators $H_e$ and $H_p$ do not satisfy the GCP.  This stems from the fact that in \cite{Cald_Ette}, (as we briefly mentioned in the introduction) the discrete analog of the gradient of a function is defined by 
		\begin{align*}
			\grad_G u(x_i) = \left(  
			u(x_i)-u(x_1),
			u(x_i)-u(x_2),
			\dots,
			u(x_i)-u(x_N)
			\right),
		\end{align*}
		and the operators there were designed to be functions, monotone increasing in the ``gradient'' variable.
		
		Thus in place of $H_e$ and $H_p$ we consider related operators, whose equations yield the same solutions, but such that these new operators also satisfy the GCP. To this end, $I_e$ and $I_p$ will be replacements for $H_e$ and $H_p$, defined by 

\begin{align}
	& I_e(u,x) \coloneqq  H_e(-u,x) = \max_{z\in G} w_{z,x}(u(z)-u(x)),\label{eqEikonal:ModifiedEikonal}\\  
	& I_p(u,x) \coloneqq  H_p(-u,x) = \frac{1}{p}\sum_{z\in G} w_{z,x} \left( (u(z)-u(x))_{+} \right)^p.\label{eqEikonal:ModifiedPEikonal}
\end{align}
Of course it is clear that $u$ solves \eqref{eqIntro:HJBGeneric} with $I=I_e$ if and only if $-u$ solves the corresponding equation with $I=H_e$ for the same $f$ and $g$ replaced by $-g$ (with a similar comment for $I_p$ and $H_p$).

The $p$-eikonal operator is interesting for many reasons, but one of them is how it demonstrates the different methods and arguments involved in the various approaches to (viscosity solutions in the continuum case) Hamilton-Jacobi equations and elliptic equations.  Our goal is to establish the following observation for $H_e$ and $H_p$ that shows they fit within our framework here.  Subsequently, we will make a small detour through some convex analysis and the arguments for the Bellman operator in the previous subsection to show that the eikonal and $p$-eikonal operators are very similar to those of Bellman.	 Before we do so, it will be useful to define a distance function given a pair $(G, w_{z,x})$.

		Given edge weights, there is a natural notion of paths and distance on $G$.  A \emph{path} is simply a function $x: \Natural\union \{0\}\to G$, and we will say $x$ is a path from $y$ to $z$, provided there exists $k\in\Natural\union \{0\}$ with $x(0)=y$ and $x(k)=z$,
		and we denote the collection of all paths from $y$ to $z$ as
		\begin{align*}
			\P(y,z) = \{ x: \Natural\union\{0\}\to G\mid x\ \text{is a path from}\ y\ \text{to}\ z\}.
		\end{align*}
		Then one can define a notion of \emph{path distance} as
		\begin{align*}
			d_{G,w_{ij}}(y,z) \coloneqq  \min_{x\in\P(y,z),\ t_z\ s.t.\ x(t_z)=z}
			\left(  \sum_{t=0}^{t_z-1} (w(x(t),x(t+1)))^{-1}  \right),
		\end{align*}
		where by convention
		\begin{align*}
			\text{if}\ \ w_{x,y}=0,\ \ \text{then}\ \ (w_{x,y})^{-1}=+\infty,
		\end{align*}
		so that the distance between two points is finite only when they can be connected by edges with non-zero weights.  This can also be used to define a notion of distance to $\Gam$,
		\begin{align*}
			d_{G,w_{ij}}(y,\Gam) \coloneqq  \min_{z\in\Gam} d_{G,w_{ij}}(y,z).
		\end{align*}
	Note however, that the path distance is not necessarily a distance in the usual sense. First, because of the above convention, it is possible for $d_{G,w_{ij}}(y,z)=\infty$ if $y$ and $z$ are contained in separate connected components of $G$. Second, it is possible that edges have nonsymmetric weights (i.e., $w_{ij}\neq w_{ji}$), in which case $d_{G,w_{ij}}$ may not be symmetric.
		
		With the above in mind, we make the following assumption.
\begin{assumption}[$G$ connected]
	For the remainder of this section, we assume that $G$ and $\{w_{x,y}\}$ forms a path connected graph, or equivalently, for each $x$, $y\in G$, $d_{G, w_{ij}}(x,y)<\infty$.
\end{assumption}

		This path distance function can be used to ensure the assumptions Theorems \ref{thmIntro:Comparison}~\eqref{assume:StrictSubSolPerturbation} and~\eqref{assume:StrictSuperSolPerturbation} hold, and for appropriate $f$ that sub/supersolutions satisfying \eqref{eqn: special sub/super} from Theorem \ref{thmIntro:Existence} can be constructed, and the condition for existence as in Remark \ref{rem: alternate existence condition} holds. We will need the following key (known) properties, see also \cite[Lemma 10]{Cald_Ette} and \cite[Theorem 12]{Cald_Ette}. For completeness we provide some of the proofs.
		
\begin{lemma}[From \cite{BungertCalderRoith-2023UnifRatesLipLearningOnGraphs} and \cite{Cald_Ette}]\label{lemEikonal:SuperSolutions}
	If $d_G$ and $d_{G^{1/p}}$ are the functions given by
	\begin{align*}
		d_{G}(x) \coloneqq  d_{G,w_{ij}}(x,\Gam)\ \ \text{and}\ \ 
		d_{G^{1/p}}(x) \coloneqq  d_{G,w_{ij}^{1/p}}(x,\Gam),
	\end{align*}
	then 
	\begin{enumerate}[(i)]
		\item \label{eqn: eikonal RHS 1} $d_G$ satisfies
		\begin{align*}
			\begin{cases}
				H_e(d_G,x) = 1\ &\text{on}\ G\setminus\Gam,\\
				d_G(x)=0\ &\text{on}\ \Gam.
			\end{cases}
		\end{align*}
		\item \label{eqn: constant bound p-eikonal}$d_{G^{1/p}}$ satisfies
		\begin{align*}
		\begin{cases}
			H_p(d_{G^{1/p}},x)\geq \frac{1}{p} &\text{on}\ \ G\setminus\Gam,\\
			d_G(x)=0\ &\text{on}\ \Gam.
		\end{cases}			
		\end{align*}
	\end{enumerate}
\end{lemma}
\begin{proof}
 The claim~\eqref{eqn: eikonal RHS 1} is proved in \cite[Lemma 3.5]{BungertCalderRoith-2023UnifRatesLipLearningOnGraphs}. For~\eqref{eqn: constant bound p-eikonal}, 
\begin{align*}
 H_p(d_{G^{1/p}},x)
 &=\sum_{y\in G}\frac{1}{p}w_{y,x}((d_{G^{1/p}}(x)-d_{G^{1/p}}(y))_+)^p\\
 &\geq \frac{1}{p}\left[\max_{y\in G}\left(w_{y,x}^{1/p}(d_{G^{1/p}}(x)-d_{G^{1/p}}(y))\right)\right]^p
 =\frac{1}{p}
\end{align*}
from \eqref{eqn: eikonal RHS 1}, noting that the penultimate expression is the operator $H_e$ defined with the graph whose weights are $w_{ij}^{1/p}$ applied to $d_{G^{1/p}}$, which is also the path distance for that graph.
\end{proof}

We now turn to verification of the assumptions in our main theorems for the operators $I_e$ and $I_p$.

\begin{proposition}\label{propEikonal:SatisfyAssumptions}	
\begin{enumerate}[(i)]
 \item \label{eqn: eikonal satisfy assumptions}The modified eikonal and p-eikonal operators, $I_e$ and $I_p$ in \eqref{eqEikonal:ModifiedEikonal} and \eqref{eqEikonal:ModifiedPEikonal} satisfy Assumptions \ref{assume:Uniqueness} and \ref{assume:Existence}.
 \item \label{eqn: eikonal comparison assumptions} Any $u\in C(G)$ such that $I_e(u, \cdot)> 0$ (resp. $I_p(u, \cdot)> 0$) on $G$ satisfies Theorem \ref{thmIntro:Comparison}~\eqref{assume:StrictSubSolPerturbation}, and any $v\in C(G)$ such that $I_e(v, \cdot)>0$ (resp. $I_p(u, \cdot)> 0$) on $G$ satisfies Theorem \ref{thmIntro:Comparison}~\eqref{assume:StrictSuperSolPerturbation} for $I_e$ (resp. $I_p$).
 \item \label{eqn: eikonal existence assumption} If $f>0$ on $G$, there exist sub and supersolutions of $I_e$ and $I_p$ which satisfy assumptions \eqref{eqn: special sub/super} in Theorem \ref{thmIntro:Existence}.
\end{enumerate}
\end{proposition}

\begin{proof}[Proof of Proposition \ref{propEikonal:SatisfyAssumptions}]
	We note that by construction, both $I_e$ and $I_p$ enjoy the GCP, and furthermore both are invariant under adding a constant to $u$ (so that Assumption \ref{assume:Uniqueness}~\eqref{assume:SubtractConstant} in fact is an equality).  Hence, Assumption \ref{assume:Uniqueness} is satisfied.  Furthermore, $I_e$ and $I_p$ are continuous functions from $C(G)$ to $C(G)$, and so Assumption \ref{assume:Existence} is immediate, proving~\eqref{eqn: eikonal satisfy assumptions}.
	
	To show \eqref{eqn: eikonal comparison assumptions}, first suppose $u\in C(G)$ and $I_e(u, \cdot)>0$ or $I_p(u, \cdot)> 0$ on $G$. As $I_e$ and $I_p$ are both positively $q$-homogeneous (with $q=1$ for $I_e$ and $q=p$ for $I_p$), we can apply Lemma \ref{lemIntro:StrictSubsolPerturbation}~\eqref{assumptionLemIntro:PositiveHomogeneous} to see Theorem \ref{thmIntro:Comparison}~\eqref{assume:StrictSubSolPerturbation} holds.
	 
	 On the other hand, if $v\in C(G)$ is such that $I_e(v, \cdot)>0$ or $I_p(v, \cdot)>0$ on $G$, first we note that as maxima of convex functions, $I_e(\cdot, x)$ and $I_p(\cdot, x)$ are both convex on $C(G)$ for any fixed $x\in G$. Then we define
	\begin{align*}
	\lambda_v&\coloneqq  \frac{1}{2}\min_{G} I_e(v, \cdot)>0,\qquad
	\lambda_{v, p}\coloneqq \frac{1}{2}\min_{G} I_p(v, \cdot)>0,\\
		\phi_{v, e} &\coloneqq  -\lambda_v d_{G},\qquad
		\phi_{v, p} \coloneqq  -\left(\frac{\lambda_{v, p}}{2\max_{G}H_p(d_{G^{1/p}},\cdot)}\right)^{\frac{1}{p}}d_{G^{1/p}};
	\end{align*}
	clearly $H_p(d_{G^{1/p}},\cdot)$ is finite valued on $G$, thus by  Lemma \ref{lemEikonal:SuperSolutions}~\eqref{eqn: constant bound p-eikonal} we see that $\max_{G}H_p(d_{G^{1/p}},\cdot)$ is positive and finite. Then by Lemma \ref{lemEikonal:SuperSolutions} and the positive $1$-homogeneity (resp. $p$-homogeneity) of $I_e$  (resp. $I_p$), 
	\begin{align*}
		I_e(\phi_{v, e}, \cdot)&=H_e(\lambda_vd_G, \cdot)=\lambda_vH_e(d_G, \cdot)= \frac{1}{2}\min_{G} I_e(v, \cdot)\leq -\lambda_v+ I_e(v, \cdot),\\
		I_p(\phi_{v, p}, \cdot)&=H_p(\lambda_{v, p}d_G^{1/p}, \cdot)=\left(\frac{\lambda_{v, p}}{2\max_{G}H_p(d_{G^{1/p}},\cdot)}\right)H_p(d_G^{1/p}, \cdot)\leq -\lambda_{v, p}+ I_p(v, \cdot).
	\end{align*}
	Thus in both cases, we may apply Remark~\ref{remIntro:ConvexSuperSolutionPerturbation} to see that $v$ satisfies Theorem \ref{thmIntro:Comparison}~\eqref{assume:StrictSuperSolPerturbation}.

	Finally we will show~\eqref{eqn: eikonal existence assumption}, using an argument similar to the one in \cite[Theorem 12]{Cald_Ette}. Since $I_e$ and $I_p$ take constant functions to the zero function, we may take $V_{f, g, \super}\equiv \min_\Gamma g$, which satisfies \eqref{eqn: special sub/super} for both operators. By using the positive $1$-homogeneity of $I_e$ and Lemma~\ref{lemEikonal:SuperSolutions}~\eqref{eqn: eikonal RHS 1}, followed by the invariance of $I_e(\cdot, x)$ under translation by constants, we can take
	\begin{align*}
		U_{f, g, \sub}\coloneqq -(\max_G f) d_{G} - \min_\Gam g
	\end{align*}
	which will satisfy \eqref{eqn: special sub/super}. By $p$-homogeneity of $I_p$, an analogous argument using Lemma~\ref{lemEikonal:SuperSolutions}~\eqref{eqn: constant bound p-eikonal} shows that we may take
	\begin{align*}
		U_{f, g, \sub}\coloneqq -(p\max_G f)^{\frac{1}{p}}d_{G^{1/p}} - \max_\Gam g.
	\end{align*}
\end{proof}

\begin{rem}
It is not clear if the operators $I_e$ and $I_p$ satisfy \eqref{eqIntro:ExistenceTheorem-ComparisonAssumption} in Theorem~\ref{thmIntro:Existence}. However, if $f>0$ the above proof shows that subsolutions with $I(u, \cdot)\geq f$ satisfy one of conditions as stated in Remark~\ref{rem: alternate existence condition}, hence one can still obtain uniqueness and existence of equations with the operator $I_e$ or $I_p$ in this case.
\end{rem}


\subsection{Some convex analysis for the $p$-eikonal operator ($p$-eikonal is a Bellman operator)}\label{subsec:ConvexForEikonal}

We would like to show that the eikonal and $p$-eikonal operators, at least when $p>1$, actually bear a strong resemblance to the Bellman operators that appeared in Subsection \ref{subsec:ControlProblem}.  To do so, we need to invoke some basic convex analysis for $I_p$.  For the remainder of this section, we assume $p\geq1$.

Let $\phi: \real^N\to\real$ be convex; again by \cite[Corollary 10.1.1]{RockafellarBook1970} such a $\phi$ is automatically continuous. Then the Legendre transform $\phi^*$ of $\phi$ and (\cite[Section 12]{RockafellarBook1970}) the subdifferential $\partial\phi(u)$ of $\phi$ at $u\in \R^N$ (\cite[Section 23]{RockafellarBook1970}) are defined by
\begin{align*}
	\phi^*(p) &\coloneqq  \sup_{u\in\real^N}\left(  u\cdot p - \phi(u)   \right),\\
	\partial \phi(u)&\coloneqq \{p\in \R^N : \phi(z)\geq \phi(u)+p\cdot (z-u),\ \forall z\in \R^N\},
\end{align*}
where we associate $\real^N$ with its dual space for the variable, $p\in\real^N$. Later, we will consider $I_p(\cdot, x): \real^N\to\real$ as a convex function for fixed $x\in G$, identifying  $u\in C(G)$ with the associated vector in $\R^N$ by an abuse of notation.

We now prove a simple lemma on local properties of the subdifferential of a convex function; below $B_R\subset \R^N$ is the open ball of radius $R$ centered at the origin.
\begin{lemma}\label{lemConvexEikonal:ConstructFromGradientsAndDual}
	
	If $\phi: \R^N\to \R$ is convex, then for each $R>0$ and any $u\in B_R$,
	\begin{align*}
		\phi(u) = \sup_{p\in\partial\phi(B_R)}\left(  p\cdot u - \phi^*(p)  \right).
	\end{align*}

\end{lemma}

\begin{proof}[Proof of Lemma~\ref{lemConvexEikonal:ConstructFromGradientsAndDual}]
	
	By \cite[Theorem 23.5]{RockafellarBook1970}, for any $p\in \partial \phi(u)$ it holds that $\phi(u) = u\cdot p - \phi^*(p)$. 
	Thus for any $u\in B_R$,
	\begin{align*}
		\phi(u)\leq \sup_{p\in\partial\phi(B_R)}\left(   p\cdot u - \phi^*(p)   \right).
	\end{align*}
	On the other hand, by using \cite[Theorem 12.2 and Corollary 7.4.2]{RockafellarBook1970} to obtain the first equality below, we know that
	\begin{align*}
		\phi(u) = \phi^{**}(u) = \sup_{p\in\real^N}\left( p\cdot u - \phi^*(p) \right)
		\geq \sup_{p\in \partial \phi(B_R)}\left(    p\cdot u- \phi^*(p)  \right).
	\end{align*}
	Hence, we have obtained our claim.
\end{proof}

We would now like to apply Lemma~\ref{lemConvexEikonal:ConstructFromGradientsAndDual} to $I_e$ and $I_p$.  We will first focus on $I_p$, and for simplicity, we will assume that $p>1$ so that for $x_i$ fixed, $I_p(\cdot,x_i):\real^N\to\real$ is differentiable. Fix $p>1$ and $x_i\in G$ and let use define $\phi_i: \real^N\to\real$ by
\begin{align*}
	\phi_i(u) \coloneqq  I_p(u,x_i),
\end{align*}
again we identify $u=(u_1, \ldots, u_N)\in \R^N$ with a function in $C(G)$ via $u_i=u(x_i)$. 
Introducing some more notation, for $i\in\{1,\dots,N\}$ let use write
\begin{align*}
	f_{ji}(u)\coloneqq \frac{1}{p}w_{ij}((u_j-u_i)_+)^p,
\end{align*}
so that
\begin{align*}
	\phi_i(u) = \sum_{j=1}^N \frac{1}{p} w_{ji}\left( (u_j-u_i)_+ \right)^p
	= \sum_{j=1}^N f_{ji}(u).
\end{align*}
Since $p>1$, each $f_{ji}$ is differentiable and a direct computation establishes that
		\begin{align*}
			\partial_kf_{ji}(u)=
			    \begin{cases}
				0&\mbox{if }k\notin\{i,j\},\\
				w_{ji}((u_j-u_i)_+)^{p-1}&\mbox{if }k=j,\\
				-w_{ji}((u_j-u_i)_+)^{p-1}&\mbox{if }k=i,
			    \end{cases}
	    \end{align*}
		which also gives 
		\begin{align}\label{eqSetting:Gradienth_i}
			\partial_k\phi_{i}(u)=
			\begin{cases}
				w_{ki}((u_k-u_i)_+)^{p-1}&\mbox{if }k\neq i,\\
				-\sum\limits_{j=1}^N w_{ji}((u_j-u_i)_+)^{p-1}&\mbox{if }k=i.
			\end{cases}
		\end{align}

		Now, as suggested by the representation in Lemma~\ref{lemConvexEikonal:ConstructFromGradientsAndDual}, if $\norm{u}\leq R$, we will want to see what happens to $\grad \phi_i(v) \cdot u$, for $\norm{v}\leq R$ (here, the variable $p$ will be of the form $p=\grad \phi_i (v)$, as $\phi_i$ is differentiable).  First, keeping in mind \cite[Lemma 2.1]{GuSc-2019MinMaxNonlocalCALCVARPDE}, we calculate for any $v$,
\begin{align*}
	\grad \phi_i(v)\cdot u &= \sum_{k=1}^N \partial_k \phi_i(v) u_k\\
	&=\sum_{k=1}^N \partial_k \phi_i(v) (u_k-u_i + u_i)\\
	&= u_i\sum_{\ell=1}^N\partial_\ell \phi_i(v)
	   +  \sum_{k=1}^N \partial_k \phi_i(v) (u_k-u_i)\\
	   &= \sum_{k\neq i}^N \partial_k \phi_i(v) (u_k-u_i)
	   =\sum_{k\neq i} w_{ki}((v_k-v_i)_+)^{p-1}(u_k-u_i).
\end{align*}
using~\eqref{eqSetting:Gradienth_i} to obtain the fourth and fifth equalities. 
By differentiability of $\phi_i$ (see \cite[Theorem 25.1]{RockafellarBook1970}) we have that $\partial \phi_i(v)=\{\nabla \phi_i(v)\}$, then combining the above calculation with Lemma~\ref{lemConvexEikonal:ConstructFromGradientsAndDual}, we see that 
\begin{align}
	\phi_i(u) &= \sup_{p\in\partial\phi_i(B_R)}\left( p\cdot u - \phi_i^*(p) \right)
	= \sup_{v\in B_R}\left( \grad\phi_i(v)\cdot u - \phi_i^*(\grad\phi_i(v)) \right)\notag\\
	&= \sup_{v\in B_R}\left( \sum_{k\not=i} w_{ki}((v_k-v_i)_+)^{p-1}(u_k-u_i) - \phi_i^*(\grad\phi_i(v)) \right).\label{eqConvexEikonal:PhiSupAlmostThere}
\end{align}

Note the similarities of the above to the Bellman type operator of Subsection \ref{subsec:ControlProblem}.  Here we only intend to point out an analogy, and so in what follows, we ignore the issue of whether or not the corresponding $K^\al$ have unit mass.  Let us introduce some more notation to make it look slightly more similar.  We want to go back from the notation of $x_i,x_j\in G$ and $u(x_i)=u_i$, and so to this end, if $v\in C(G)$, we will write
\begin{align*}
	v_i\coloneqq v(x_i)\ \ \ \text{and}\ \ \ v_x\coloneqq v(x),
\end{align*}
and the particular choice should be clear from context.  Furthermore, we let
\begin{align*}
	& \A \subset C(G),\ \ \A \coloneqq  B_R = \{ v\in C(G)\mid \norm{v}\leq R \},\\
	& K_v(x,y) \coloneqq  w(y,x)\left( (v(y)-v(x))_+ \right)^{p-1},\\
	& f_v(x) \coloneqq  \phi_i^*(\grad\phi_i(v))\ \text{when}\ x=x_i.
\end{align*}
Thus, we see that the representation in \eqref{eqConvexEikonal:PhiSupAlmostThere} gives, at least for $p>1$,
\begin{align}\label{eqConvexEikonal:PEikonalIsSupOfLinear}
	\text{if}\ u\in B_R,\ \ 
	I_p(u,x) &= 
	\sup_{v\in B_R}\left( -f_v(x)  +  \sum_{y\in G} w(y,x)\left( (v(y)-v(x))_+ \right)^{p-1}(u(y)-u(x))    \right) \nonumber      \\
	&= \sup_{v\in B_R}\left( -f_v(x) +  \sum_{y\in G} K_v(y,x)(u(y)-u(x))  \right).
\end{align}
We note that the presence of $f_v$ inside of the supremum can be incorporated into the structure from subsection \ref{subsec:ControlProblem}, whereby the running cost, $f$, would also be allowed to depend upon $\al$.   Thus, $I_p$ has the structure of a Bellman operator.  It should be noted that if one were to put this directly into the context of subsection \ref{subsec:ControlProblem} that the control problem is very degenerate.

We can make a similar representation for $I_e$.  The parameters in this setting are as follows:
\begin{align*}
	\A \coloneqq  G,
\end{align*}
and
\begin{align*}
	I_e(u,x) = \max_{z\in G} \sum_{y\in G}K_z(x,y)(u(y)-u(x))
	= \max_{z\in G} w_{z,x}(u(z)-u(x)),
\end{align*}
so that
\begin{align*}
	K_z(x,y) \coloneqq  
	\begin{cases}
		w_{y,x}\ &\text{if}\ y=z,\\
		0\ &\text{if}\ y\not=z.
	\end{cases}
\end{align*}

One can see, by inspection, that the collection of linear operators that appear in the formulation of $I_p$ is much richer than those allowed in the representation of $I_e$.  Hence, it would be expected that the operator, $I_p$, would be less sensitive to changes in the weights, $w_{x,y}$, than is the operator, $I_e$.  This robustness was already observed by different means in \cite{Cald_Ette}, and indeed, this was one of the most important features of $I_p$.

		\subsection{Operators similar to the $p$-eikonal operator}

		Motivated by the structure of the $p$-eikonal operator, we also consider a related family, which in some settings includes the $p$-eikonal operator as a particular example.  We keep the edge weight structure as given in the previous subsection, as in \eqref{eqEikonal:EdgeWeightNotation}.  Then a class of operators with similar structure to the $p$-eikonal operator are of the form
		\begin{align*}
			J(u,x)\coloneqq  \sum_{y\in G}w_{y,x}c\left( u(y)-u(x) \right),
		\end{align*}
		where $c:\real\to\real$ is continuous and increasing.  We note that similar nonlocal operators in the continuum setting of integro-differential equations were already suggested as examples of the fully nonlinear theory in \cite[Sections 3 and 6]{CaSi-09RegularityIntegroDiff}.  Unless $c$ has special structure (such as convexity), $J$ will not be a convex operator.  Instead, $J$ is typically a good example of a fully nonlinear elliptic operator on $C(G)$.

		Under just the assumption that $c$ is continuous and increasing, $J$ will be continuous on $C(G)$, obeys the GCP, and is invariant with respect to addition of constants.  Thus, Assumptions \ref{assume:Uniqueness} and \ref{assume:Existence} are easily confirmed to be satisfied.  The much trickier part is to know whether or not the additional assumptions included in Theorems \ref{thmIntro:Comparison} and \ref{thmIntro:Existence} are also verified (the perturbation to a strict sub/supersolution, and the existence of a particular pair of subsolution and supersolution that can be used in the Perron construction).  In the case of the $p$-eikonal operator, two factors played a role in verifying these additional assumptions, and they came from the very particular structure of $I_p$; namely how the path distance function behaves with respect to $I_p$, and the fact that $I_p$ is convex.  Indeed, one will note that even though we demonstrated that $I_p$ has the structure of a Bellman operator as treated in subsection \ref{subsec:ControlProblem}, we took advantage of the convexity of $I_p$ and the path distance function to confirm the assumptions required for existence and uniqueness, instead of relying on the arguments for the Bellman operator in subsection \ref{subsec:ControlProblem}.
		
		In this setting, thought of as a nonlocal elliptic operator, if $c$ is differentiable on $\R$ the values of $c'\geq 0$ somehow dictate how ``degenerate'' the operator, $J$, will be.  That is to say, if it is allowed that $c'(t)=0$ for some $t\in\real$, then this is rather degenerate, and it is harder to verify the additional assumptions included in Theorems \ref{thmIntro:Comparison} and \ref{thmIntro:Existence}.  For example, in the setting of $I_p$, the corresponding $c$ satisfies $c'(t)=0$ for $t\leq 0$.  In contrast,  one can add the following ``uniform ellipticity'' assumption on $c$, given by requiring $c$ is Lipschitz in $t$ and
\begin{align*}
	\exists\ \lam, \Lam>0,\ \ \forall\ t\in\real,\ \ \lam\leq c'(t)\leq \Lam. 
\end{align*}
Note that the $c$ corresponding to $I_p$ does \emph{not} satisfy this condition.

In the uniformly elliptic setting of nonlocal equations, the approach to verifying the perturbation of subsolutions and the existence of particular sub and supersolutions often follows along the lines mentioned in subsection \ref{subsec:ControlProblem} with regards to the Bellman operator, the property given in \eqref{eqControl:MinimalMaximalInequalities}, and the approach to proving Propositions \ref{propControl:SatisfiesAssumptionUniqueness} and \ref{propControl:SatisfiesAssumptionExistence}.  To this end, if we consider a generic function, $\phi\in C(G)$, by the mean value theorem there exists $z\in \real$, with 
\begin{align*}
	z = z(x,y,u,\phi),\ \ \text{and}\ \ 
	J(u + \phi,x) - J(u,x) = 
	\sum_{y\in G} w_{y,x}c'(z)(\phi(y)-\phi(x)).
\end{align*}
Not surprisingly, if one thinks of the choice of this $z$ as some black box operation, then at least formally, this operator on $\phi$ is quite similar to the graph Laplacian appearing in subsection \ref{subsec:Markov}.  Making some further manipulations, we see that
\begin{align*}
	J(u + \phi,x) - J(u,x) = 
	\sum_{y\in G} w_{y,x}c'(z)(\phi(y)-\phi(x))
	\geq \inf_{\lam\leq a\leq \Lam} \sum_{y\in G} a\cdot w_{y,x}(\phi(y)-\phi(x)).
\end{align*}
Since we have allowed for any $a\in[\lam,\Lam]$, one can check that
\begin{align}\label{eqGeneralEikonal:DefMinimalOpForJ}
	\M_J^-(\phi,x)&\coloneqq 
	\inf_{\lam\leq a\leq \Lam} \sum_{y\in G} a\cdot w_{y,x}(\phi(y)-\phi(x))\\
	&= \sum_{y\in G} \lam w_{y,x}(\phi(y)-\phi(x))_+ 
	+ \Lam w_{y,x}(\phi(y)-\phi(x))_{-},\nonumber
\end{align}
which again, bears a very strong similarity to the example of the minimal operator appearing in \cite[Section 3]{CaSi-09RegularityIntegroDiff}.  Up to a normalization constant, this also has the same structure as the Bellman operator in subsection \ref{subsec:ControlProblem}.  Indeed, if one assumes that the edge weights $w_{y,x}$ have unit mass, the question of satisfying the additional assumptions of Theorems \ref{thmIntro:Comparison} and \ref{thmIntro:Existence} comes down to verifying whether or not the exit times in Assumption \ref{assumeControl:FiniteExitTime} are finite, which is equivalent to solving the minimal equation by Lemma~\ref{lemControl:FiniteExitEquivSolBellmanRHS1}.  In this case, this will reduce to checking the structure of the Markov chain induced by $w_{y,x}$, as in subsection \ref{subsec:Markov}.

We can see, by inspecting the structure of the linear operators appearing in each of $I_e$, $I_p$, and the minimal operator associated to $J$ (i.e. $\M_J^-$ in \eqref{eqGeneralEikonal:DefMinimalOpForJ}), that these each have a distinct difference in their level of degeneracy, or rather, their sensitivity with respect to changes in $w_{x,y}$.

		\subsection{General min-max}\label{subsec:GeneralMinMax}

		The structure encountered in the Bellman operator \eqref{eqControl:BellmanOperatorInf} and the $p$-eikonal operator \eqref{eqConvexEikonal:PEikonalIsSupOfLinear} has a natural generalization to non-convex (non-concave) settings, and those are operators that take the form:
		\begin{align}\label{eqGeneralMinMax:FirstMinMax}
			I(u,x)=\min_{\al\in\A}\max_{\beta\in\B_\al}\left\{f^{\al\beta}(x)+c^{\al\beta}(x)u(x)+\sum_{y\in G}K^{\al\beta}(x,y)(u(y)-u(x))\right\},
		\end{align}
		where $\A$ is an index set, and for each $\al$, $\B_\al$ is an index set which may depend on $\al$, and for each $\al,\beta$, $f^{\al\beta}, c^{\al\beta}\in C(G)$, $K^{\al\beta}:G\times G\to[0,\infty)$.   Note we have included the assumption $K^{\al\beta}\geq0$, but, in fact without assuming this \emph{a priori}, it is mostly straightforward to check that if $I$ has the min-max form in \eqref{eqGeneralMinMax:FirstMinMax} and also has the GCP, then necessarily $K^{\al\beta}\geq 0$.  The operators appearing in (\ref{eqGeneralMinMax:FirstMinMax}) are often times associated with two-player games, and in that context they are called Isaacs equations.

		In fact, when $I:C(G)\to C(G)$ happens to be \emph{locally} Lipschitz on $C(G)$, then the following result, below, was observed in \cite{GuSc-2019MinMaxNonlocalCALCVARPDE} and \cite{GuSc-2019MinMaxEuclideanNATMA}, and it shows that \eqref{eqGeneralMinMax:FirstMinMax} is in fact generic for locally Lipschitz operators.  We note that all of the examples given above will be locally Lipschitz (for the Bellman operator in subsection \ref{subsec:ControlProblem}, there need to be natural assumptions, which are related to Assumption \ref{assumeControl:FiniteExitTime}).  To state the result, we mention some notation for Lipschitz functions on Banach spaces, which is similar to the notion of subdifferential for convex functions  which appeared in subsection \ref{subsec:ConvexForEikonal} (see \cite[Section 12]{RockafellarBook1970}).  Here we rely on the notion of the generalized gradient from Clarke (e.g. \cite[Chapter 2]{Clarke-1990OptimizationNonsmoothAnalysisSIAMreprint}).  If $X$ is a Banach space, and $F:X\to\real$ is Lipschitz, then the following limit defines a notion of upper gradient,
		\begin{align*}
			F^\circ(x;v) \coloneqq  \limsup_{y\to x}\limsup_{t\to 0^+}
			\frac{F(y+tv)-F(y)}{t},
		\end{align*}
		and the generalized gradient at $x$ is the set of supporting linear functionals at $x$,
		\begin{align*}
			\partial F(x) \coloneqq  \left\{  \ell\in X^*\mid F^\circ(x;v)\geq \ell(v)  \right\}.
		\end{align*}
		We let $\partial F$ be shorthand for the image of the gradient mapping, i.e. if $Y\subset X$,
		\begin{align*}
			\partial F(Y) \coloneqq  \overline{\text{c.h.}} \Union_{y\in Y}\partial F(y)\subset X^*,
		\end{align*}
		where $\overline{\text{c.h.}}$ stands for the closed convex hull.
With a few observations involving a generalization of the Mean Value Theorem to such $F$ (known as Lebourg's theorem \cite[Theorem 2.3.7]{Clarke-1990OptimizationNonsmoothAnalysisSIAMreprint}), one can deduce the following result as in \cite[Section 2]{GuSc-2019MinMaxEuclideanNATMA}
\begin{align*}
	&\text{if}\ Y\ \text{is a bounded convex subset of}\ X,\ \text{then}\ \forall\ y\in Y,\\
	&\ \ F(y) = \min_{z\in Y}\sup_{\ell\in \partial F(Y)} (\ell(y-z) + F(z)).
\end{align*}
(In fact, the supremum is always a maximum in this situation, but that is irrelevant for our usage, as we are in a finite dimensional space.)
Even though this discussion has assumed $F$ is real valued, since $C(G)$ is a finite dimensional space, all of this discussion readily extends to $I: C(G)\to C(G)$, and the following result is presented in \cite[Section 2]{GuSc-2019MinMaxNonlocalCALCVARPDE}.  It relies on the fact that
\begin{align*}
	& L: C(G)\to C(G)\ \text{is linear and bounded}\ \text{if and only if}\\ 
	& \exists\  c_{L}\in C(G)\ \text{and}\ K_{L}: G\times G\to\real\ \text{with}\
	L(v,x) = c_{L}(x)v(x) + \sum_{y\in G} K_L(x,y)(v(y)-v(x)).
\end{align*}
Here we abuse notation and we will use $\partial I(Y)$ to mean the collection of generalized derivatives, as above. Note below, since $C(G)$ is finite dimensional and $I$ is locally Lipschitz, the inner supremum is indeed a maximum.

\begin{proposition}[min-max locally Lipschitz operators-- Lemma 2.4 of \cite{GuSc-2019MinMaxNonlocalCALCVARPDE}]\label{propGeneralMinMax:MinMaxRepresentation}
	If $I: C(G)\to C(G)$ is locally Lipschitz and $Y$ is a bounded convex subset of $C(G)$, then 
	\begin{align*}
		\forall\ u\in Y,\ \ 
		I(u,x) = \min_{v\in Y}\max_{L\in \partial I(Y)} 
		\left( f_{L,v}(x) +  c_L(x)u(x) + \sum_{y\in G} K_L(x,y)(u(y)-u(x))  \right),
	\end{align*}
	where for $L\in\partial I(Y)$, we denote $f_{L,v}(x) = I(v,x)-L(v,x)$.
	Furthermore, if $I$ has the GCP, then for all $L\in\partial I(Y)$, $K_L(x,y)\geq 0$.  Finally, if for all constants $c\geq 0$ and $u\in C(G)$, $I$ satisfies $I(u+c,x)\leq I(u,x)$, then for each $L\in \partial I(B_R)$, we have $c_L\leq 0$.
\end{proposition}

\begin{rem}\label{remGeneralMinMax:KxxZero}
	We make an important note that just as was done explicitly for the kernels, $K_v$ arising at the end of subsection \ref{subsec:ConvexForEikonal}, the construction given in \cite{GuSc-2019MinMaxNonlocalCALCVARPDE} to produce the above representation necessarily gives
	\begin{align*}
		\forall\ L\in\partial I(Y),\ \forall\ x\in G,\ K_L(x,x)=0.
	\end{align*}
\end{rem}

		\subsection{A generic operator $H$ from \cite{Cald_Ette}}\label{subsec:GeneralH-CE}

		In this section we want to demonstrate the correspondence between the Hamilton-Jacobi operators considered in \cite[Section 2]{Cald_Ette} and the operators treated in our work.  This has two parts.  First, we assume that a Hamilton-Jacobi operator, $H$ is given, without any assumptions about how it satisfies the assumptions of \cite[Section 2]{Cald_Ette}.  This $H$ allows for a natural construction of $I$ which fits into our framework, which is below, in (\ref{eqGeneralH-CE:GivenHDefineI}).  We show that $H$ and $I$ equivalently satisfy the assumptions of their respective settings, Proposition \ref{propGeneralH-CE:GivenHDefineIEquivAssumptions}.  However, there is another interesting result in this direction, and that is when some $I$ is given, how does it define an $H$ that will be compatible with \cite[Section 2]{Cald_Ette}?  This is the content of Proposition \ref{propGeneralH-CE:GivenIDefH}.  The conclusion is that at least in the case of locally Lipschitz operators on $C(G)$, all of those which have a comparison property must have a min-max representation, and so any locally Lipschitz $I$ will also give a corresponding $H$.

As mentioned above, a general form of Hamilton-Jacobi equation from \cite{Cald_Ette} is given by an operator, $H$, with
\begin{align*}
	H: \real^N\times \real\times G \to \real,
\end{align*}
and $H$ is continuous and monotone in its first two variables.  The ``gradient'' of $u\in C(G)$ is given, for $x\in G$ and $G=\{x_1,\dots, x_N\}$, by 
\begin{align*}
	\grad_G u(x) = (
	u(x)-u(x_1),
	u(x)-u(x_2),
	\dots,
	u(x)-u(x_N)
	).
\end{align*}
The equation can then be written as
\begin{align}\label{eqGeneralH-CE:CalderHJB}
	\begin{cases}
		H(\grad_G u(x), u(x), x) = f(x)\ &\text{in}\ G\setminus\Gam,\\
		u=g\ &\text{on}\ \Gam.
	\end{cases}
\end{align}
Given such an $H$, we now define $I: C(G)\to C(G)$ to fit within our context by
\begin{align}\label{eqGeneralH-CE:GivenHDefineI}
	I(u,x) \coloneqq H(-\grad_G u(x), -u(x), x).
\end{align}

We note here, operators $H$ in \cite{Cald_Ette} are assumed to be monotone in the gradient variable, which is a natural assumption, however the proofs of \cite{Cald_Ette} do not actually use the full strength of this monotonicity. Rather, a very particular notion of monotonicity is used, one that applies with vectors whose $i$-th entry is zero.  That is to say the relevant property is the following.

\begin{definition}\label{defGeneralH-CE:DifferenceMonotone}
	For $p\in\real^N$, we say that $p$ is a difference-gradient at $x_i\in G$ provided $p_i=0$.
	A function, $H:\real^N\times\real\times G\to\real$ is said to be ``differences monotone'' provided that for all $x_i\in G$, if $p$ and $q$ are difference-gradients at $x_i$ with $p\leq q$, and $s\leq t$, then $H(p,s,x_i)\leq H(q,t,x_i)$.
\end{definition}

The reason for this definition is that we are studying equations that have the form (\ref{eqGeneralH-CE:CalderHJB}), and so the function $H$ is only used at $x_i$ for $p$ who have $p_i=0$ (as the difference based at $x_i$ gives a zero entry in $\grad_G u(x_i)$).  As mentioned, we hope to make an equivalence between the operators, $I$, in this work, and those $H$ in \cite{Cald_Ette}.  Since our operators use the GCP in place of monotonicity, we cannot obtain monotonicity of the function $I$ for a generic ordered pair, $p$, $q\in\real^N$.  Rather, we can only obtain information when $p$ and $q$ arise as gradients of some function, which necessarily satisfies
\begin{align*}
	\forall\ u\in C(G),\ \forall\ x_i\in G,\ (\grad_G u(x_i))_i = 0,
\end{align*}  
i.e. when evaluated at $x_i$, the $i$-th entry of $\grad_G u(x_i)$ is always zero.

This does not change the existence and uniqueness theory from \cite{Cald_Ette}, and we note the following remark without proof.

\begin{rem}\label{remGeneralH-CE:ReplaceMonotoneByDifferencesMonotone}
	All of the results in \cite[Section 2]{Cald_Ette} remain valid if all requirements that $H$ be monotone are replaced by the requirement that $H$ is ``differences monotone'' as given in Definition \ref{defGeneralH-CE:DifferenceMonotone}.
\end{rem}

It is now immediate to verify the following equivalence of assumptions for such an $I$, given an $H$, between this work and \cite{Cald_Ette}.

\begin{proposition}\label{propGeneralH-CE:GivenHDefineIEquivAssumptions}
	Assume $H:\real^N\times\real\times G\to\real$ is given and $I$ is defined by \eqref{eqGeneralH-CE:GivenHDefineI}.   Then $I$ satisfies our assumptions of uniqueness and existence if and only if $H$ satisfies the assumptions of uniqueness and existence in \cite[Theorem 7]{Cald_Ette}, with ``monotone'' replaced by ``differences monotone''.
\end{proposition}

\begin{rem}
	Note the way the previous proposition is stated, we are only assuming that $H$ is a function on the given domain.  The properties of continuity and differences monotone will be deduced from the properties of $I$ given in this work.
\end{rem}

\begin{proof}[Proof of Proposition \ref{propGeneralH-CE:GivenHDefineIEquivAssumptions}]
	
This was already worked out for the special cases of the eikonal and $p$-eikonal operators in subsection \ref{subsec:Eikonal}.  We do not work through all of the details here, but we make some notes.

First of all, assume that $I$ satisfies Assumptions \ref{assume:Uniqueness} and \ref{assume:Existence} herein.  By Lemma \ref{lemIntro:IContinuous}, $I$ will be continuous in $u$, hence $H$ will be continuous in its first two variables.  We still need to establish the differences monotonicity property of $H$.   To this end, let $x_i\in G$ be fixed, let $p$, $q\in\real^N$ be difference-gradients at $x_i$ (see Definition \ref{defGeneralH-CE:DifferenceMonotone}), and let $s$, $t\in\real$, with
\begin{align}\label{eqn: ordered}
	p\leq q\ \ \text{and}\ \ s\leq t.
\end{align} 

We now define $u$, $v\in C(G)$ by
\begin{align*}
u(x_j)=-s+p_j,\quad
v(x_j)=-t+q_j,
\end{align*} 
for all $j=1, \ldots, N$, then since $p_i=q_i=0$ we have
\begin{align*}
	\grad_G u(x_i) = -p\ \ \text{and}\ \ u(x_i)=-s,\\
	\grad_G v(x_i) = -q\ \ \text{and}\ \ v(x_i)=-t.
\end{align*}
We can then confirm by inspection that defining the function $w$ by
\begin{align*}
	w \coloneqq u - (t-s),
\end{align*}
we have
\begin{align*}
	\grad_G w(x_i) = \grad_G u(x_i) \geq \grad_G v(x_i)\ \ \text{and}\ \ 
	w(x_i) = v(x_i).
\end{align*}
These properties then give that
\begin{align*}
	\forall\ z\in G,\ w(z)\leq v(z)\ \ \text{and}\ \ w(x_i)=v(x_i),
\end{align*}
and since we assume $I$ has the GCP,
\begin{align*}
	I(w,x_i)\leq I(v,x_i).
\end{align*}
However, we see that $t-s\geq 0$ from~\eqref{eqn: ordered}, so by Assumption \ref{assume:Uniqueness}~\eqref{assume:SubtractConstant},
\begin{align*}
	I(u,x_i)\leq I(u-(t-s),x_i) = I(w,x_i).
\end{align*}
Since we also have $\grad_G w = \grad_G u$, we can combine all of these observations to see that
\begin{align*}
	H(p,s,x_i) &= H(-\grad_G u(x_i), -u(x_i),x_i)
	= I(u,x_i)\leq I(w,x_i)\\
	&\leq I(v,x_i)
	= H(-\grad_G v(x_i), -v(x_i),x_i)
	= H(q,t,x_i).
\end{align*}

Going the other direction, we only need confirm that $I$ satisfies the GCP.  However, this is immediate from the fact that 
\begin{align*}
	\text{if}\ u\leq v\ \text{and}\ u(x)=v(x),\ \text{then}
	-\grad_G u(x)\leq -\grad_G v(x),
\end{align*}
combined with the differences monotonicity of $H$.

\end{proof}

A slightly more interesting result comes from Proposition \ref{propGeneralMinMax:MinMaxRepresentation}, whereby we take a reasonably generic operator, $I$, and show that the corresponding $H$ must have the form as in \cite{Cald_Ette}.

\begin{proposition}\label{propGeneralH-CE:GivenIDefH}
	Assume that $I$ is locally Lipschitz on $C(G)$, that $f_{L,v}$, $c_L$, and $K_L$ are as in Proposition \ref{propGeneralMinMax:MinMaxRepresentation} (with $K_L(x,x)=0$), and that for each $R>0$, the function $H_R: \real^N\times \real \times G\to \real$ is defined, using the given $I$, through its min-max representation by
	\begin{align*}
		H_R(p,s,x) = 
		\min_{v\in B_R}\max_{L\in\partial I(B_R)}
		\left( f_{L,v}(x) - c_L(x)s + \sum_{x_i\in G}K_L(x,x_i)p_i \right).
	\end{align*}
	Then for each $x$ fixed, $H_R(\cdot,\cdot,x)$ is a Lipschitz function on $\real^N\times\real$.
	If, additionally, $I$ satisfies Assumption \ref{assume:Uniqueness}, then for each $x$, $H_R(p,s,x)$ is differences monotone in $p,s$ (see Definition \ref{defGeneralH-CE:DifferenceMonotone}).  Finally, if $\gam:[0,\infty)\to[0,\infty)$ is non-decreasing with $\gam(0)=0$, and for all $x\in G$, $u\in C(G)$, and $c\in [0, \infty)$,
	\begin{align*}
		I(u-c,x) - I(u,x)\geq \gam(c),
	\end{align*}
	then, for each $p\in\real^N$ and $s\leq t$, it holds that
	\begin{align*}
		H_R(p,t,x)-H_R(p,s,x)\geq \gam(t-s),
	\end{align*}
	i.e. $H_R$ is proper in the sense of \cite[Definition 5]{Cald_Ette}.
\end{proposition}

\begin{rem}
	The local Lipschitz assumption on $I$ means that that the family of ingredients in the min-max for $H_R$ naturally depends upon $R$ in the fact that the min and max are only taken only over $v\in B_R$ and $L\in\partial I(B_R)$, as in Proposition \ref{propGeneralMinMax:MinMaxRepresentation}.  However, as one can see by Proposition \ref{propGeneralMinMax:MinMaxRepresentation}, whenever $R_1<R_2$ and $(p,s,x)\in B_{R_1}\times[-R_1,R_1]\times G$, we see that
	\begin{align*}
		H_{R_2}(p,s,x) = H_{R_1}(p,s,x).
	\end{align*}
	This is like approximating the function, $f:\real^d\to\real$ with $f(x)=\frac{1}{2}\abs{x}^2$ by restricting its Legendre transform to a bounded set.  In this context, one has
	\begin{align*}
		f(x) = \sup_{y\in\real^d}\left( y\cdot x  - \frac{1}{2}\abs{x-y}^2 \right)\ \ \ \text{and}\ \ \
		f_R(x) = \sup_{\abs{y}\leq R}\left( y\cdot x  - \frac{1}{2}\abs{x-y}^2 \right).
	\end{align*}
	One can see that $f$ is only locally Lipschitz, whereas $f_R$ is globally Lipschitz, and that the two functions will agree whenever $\abs{x}\leq R$.
\end{rem}

Before we give a sketch for the proof of Proposition \ref{propGeneralH-CE:GivenIDefH},
we note the following connection between $H$ and $I$, which confirms what was to be expected from reversing the direction in the definition given by \eqref{eqGeneralH-CE:GivenHDefineI}. 

\begin{proposition}\label{propGeneralH-CE:ConfirmHEqualsI}
	 For $x\in G$, $p\in\real^N$, and $s\in\real$ fixed, if we define $v\in C(G)$ by  
	\begin{align*}	
v(y)\coloneqq
\begin{cases}
s,&y=x,\\
s-p_j,&y=x_j\neq x, 
\end{cases}
	\end{align*}
	and $R>0$ is such that $v\in B_R\subset C(G)$, then for $H_R$  defined in Proposition \ref{propGeneralH-CE:GivenIDefH},
	\begin{align*}
		H_R(p,s,x) = I(-v,x).
	\end{align*}
	
\end{proposition}
\begin{proof}[Proof of Proposition \ref{propGeneralH-CE:ConfirmHEqualsI}]
	
	First, we make an important remark about the type of $H_R$ that are obtained by this construction.  Because of the very specific structure given in the min-max, with $K_L(x,x)=0$ (see Remark \ref{remGeneralMinMax:KxxZero}), we see that at $x_i$, whenever $p,q\in\real^N$ with $p_j=q_j$ for $j\not=i$, we have that
	\begin{align*}
		H_R(p,s,x_i)=H_R(q,s,x_i)
	\end{align*}

	We note that for $x_j\neq x$, we have $(\grad_G v(x))_j = p_j$, and for $R$ large enough with $v\in B_R$,
	\begin{align*}
		H_R(p,s,x)&=
		H_R(\grad_G v(x), v(x), x) \\
		&= \min_{w\in B_R}\max_{L\in\partial I(B_R)}
		\left( f_{L,w}(x) - c_L(x)s + \sum_{x_i\in G}K_L(x,y)p_i \right)\\
		&= \min_{w\in B_R}\max_{L\in\partial I(B_R)}
		\left( f_{L,w}(x) + c_L(x)(-v(x)) + \sum_{x_i\in G}K_L(x,y)(-(v(y)-v(x))) \right)\\
		&= I(-v,x),
	\end{align*}
	where the last equality is from Proposition \ref{propGeneralMinMax:MinMaxRepresentation}.

\end{proof}

Now we will give a sketch of the argument for Proposition \ref{propGeneralH-CE:GivenIDefH}.

\begin{proof}[Sketch of the proof of Proposition \ref{propGeneralH-CE:GivenIDefH}]

	First, we show that $H_R(\cdot,\cdot,x)$ is Lipschitz on $B_R\subseteq \real^N\times\real$ for a fixed $x\in G$.	
	Indeed, for $x\in G$, $p$, $q\in B_R$, and $s$, $t\in [-R, R]$, by the definition of $H$ we find
	
\begin{align*}
 H(p, s, x)-H(q, t, x)
 &\leq \min_{v_1\in B_R}\max_{L_1\in\partial I(B_R)}\max_{v_2\in B_R}\min_{L_2\in\partial I(B_R)}
		\left(f_{L_1, v_1}(x)- c_{L_1}(x)s + \sum_{x_i\in G}K_{L_1}(x,x_i)p_i\right.\\
		&\left. -f_{L_2, v_2}(x)+ c_{L_2}(x)t - \sum_{x_i\in G}K_{L_2}(x,x_i)q_i\right)\\
		&\leq \max_{L_1\in\partial I(B_R)}\max_{v_2\in B_R}
		\left(c_{L_1}(x)(s-t) + \sum_{x_i\in G}K_{L_1}(x,x_i)(p_i-q_i)\right)\\
		&\leq \max_{L_1\in\partial I(B_R)}
		\left(\lvert c_{L_1}(x)\rvert + \sum_{x_i\in G}\lvert K_{L_1}(x,x_i)\rvert\right) (\lvert s-t\rvert +\lvert p-q\rvert).
\end{align*}
To conclude, we note that by construction, see \cite[Section 2]{GuSc-2019MinMaxNonlocalCALCVARPDE}, since $I$ is Lipschitz on $B_R$, all operators in $\partial I(B_R)$ have an operator norm bounded by the Lipschitz constant of $I$ on $B_R$.  Thus, $H_R$ is Lipschitz.

	Next we discuss the monotonicity of $H$.  With regards to the $p$ variable, this comes from the fact that $K_L\geq 0$ (which is a consequence of the GCP for $I$, see \cite[Section 2]{GuSc-2019MinMaxNonlocalCALCVARPDE}).
	The monotonicity in $s$ is due to the fact that $c_L\leq 0$, which is a result of Assumption \ref{assume:SubtractConstant} and the fact that $\partial I(B_R)$ is a closed convex hull of operators obtained as derivatives of $I$ (see \cite[Section 2]{GuSc-2019MinMaxNonlocalCALCVARPDE}).
	 So, we note the following at all $u\in C(G)$ for which $I$ happens to be differentiable (which is, by Rademacher's Theorem, a.e. $u$). We note that after reversing the subtraction of a constant in Assumption \ref{assume:Uniqueness}, we have for $c\geq 0$,
	\begin{align*}
		\text{if}\ t>0,\ \text{then}\ 
		(I(u+tc,x) - I(u,x))\leq 0
		\ \ \ \text{and}\ \ \ 
		\text{if}\ t<0,\ \text{then}\ 
		(I(u+tc,x) - I(u,x))\geq 0.
	\end{align*}
Thus,  if $L$ happens to be a derivative of $I$ at $u$,
	\begin{align*}
		L(c,x) = \lim_{t\to0}\frac{1}{t}(I(u+tc,x) - I(u,x))\leq 0.
	\end{align*}
	But
	\begin{align*}
		L(c,x) = c_L(x)\cdot c\leq 0.
	\end{align*}
	So $c_L\leq 0$. (Note, we only care about this observation for $c=1$.)

	Finally, we address the proper condition. 
	Suppose there is a non-decreasing function $\gam\geq 0$, $\gam(0)=0$, so that for all $x\in G$ and $c \in [0, \infty)$
	\begin{align*}
		I(u-c,x) - I(u,x)\geq \gam(c).
	\end{align*}
	We let $p\in\real^N$ be fixed, and $s\leq t$, and we then define $u$, $v\in C(G)$ as
	\begin{align*}
		u(x) = s,\ \ u(x_i)= s - p_i
		\ \ \ \text{and}\ \ \
		v(x) = t,\ \ v(x_i) = t - p_i.
	\end{align*}
	Then since $t-s\geq 0$ and $u = v -(t-s)$, using Proposition \ref{propGeneralH-CE:ConfirmHEqualsI} we obtain
	\begin{align*}
		H(p,t,x) = I(-v,x) = I(-u - (t-s), x) \geq I(-u,x) + \gam(t-s) = H(p,s,x) + \gam(t-s).
	\end{align*}

\end{proof}



	\bibliography{./refs-HJB-graph.bib}
	\bibliographystyle{alpha}

\end{document}